\documentclass[12pt]{article}
\usepackage{amsmath}
\usepackage{amssymb}
\numberwithin{equation}{section}
\usepackage{epic}

\title{Joint probability for the
 Pearcey process}

\author{Mark Adler\thanks{2000
{\em Mathematics Subject Classification}. Primary:
60G60, 60G65, 35Q53; secondary: 60G10, 35Q58. {\em Key
words and Phrases}: Dyson's Brownian motion, Pearcey
process, random matrices coupled in a chain, random
matrices with external potential, infinite-dimensional
diffusions.
\newline
 Department of Mathematics, Brandeis University,
Waltham, Mass 02454, USA. E-mail: adler@brandeis.edu.
The support of a National Science Foundation grant \#
DMS-01-00782 is gratefully acknowledged.}~~~~~~ Pierre
van Moerbeke\thanks{ Department of Mathematics,
Universit\'e de Louvain, 1348 Louvain-la-Neuve, Belgium
and Brandeis University, Waltham, Mass 02454, USA. This
work was done while PvM was a Miller visiting Professor
at the University of California, Berkeley, USA E-mail:
vanmoerbeke@math.ucl.ac.be and @brandeis.edu. The
support of a National Science Foundation grant \#
DMS-04-06287, a European Science Foundation grant
(MISGAM), a Marie Curie Grant (ENIGMA), FNRS and
Francqui Foundation grants is gratefully
acknowledged..}}

\date{}

\newcommand{\MAT}[1]{\left(\begin{array}{*#1c}}
\newcommand{\mat}{\end{array}\right)}
\newcommand{\qed}{\leavevmode\unskip\nobreak\penalty200\hskip2pt\null
\nobreak\hfill\rule{1.1ex}{1.1ex}
\medbreak }

\newcommand{\rg}{\rightarrow}

\newcommand{\AR}{{\cal A}}
\newcommand{\CR}{{\cal C}}
\newcommand{\DR}{{\cal D}}

\newcommand{\LR}{{\cal L}}
\newcommand{\MR}{{\cal M}}

\newcommand{\PR}{{\cal P}}

\newcommand{\BB}{{\cal B}}
\newcommand{\BT}{{\cal T}}
\newcommand{\XR}{{\cal X}}
\newcommand{\BC}{{\mathbb C}}

\newcommand{\BP}{{\mathbb P}}
\newcommand{\BQ}{{\mathbb Q}}

\newcommand{\iy}{\infty}
\newcommand{\pl}{\partial}
\newcommand{\al}{\alpha}

\newcommand{\no}{\nonumber}

\def\inn#1#2{\left\langle#1\,\left\vert\,#2\right.\right\rangle}

\newenvironment{proof}{\medskip\noindent{\it Proof:\/} }{\qed}
\newenvironment
        {remark}{\medskip\noindent\underline{\it Remark:\/} }{\medbreak}

\newcommand{\la}{\langle}
\newcommand{\ra}{\rangle}
\newcommand{\ga}{\gamma}

\newcommand{\dt}{\delta}
\newcommand{\Dt}{\Delta}
 \newcommand{\vr}{\varepsilon}
\newcommand{\sg}{\sigma}
\newcommand{\BR}{{\mathbb R}}

\newcommand{\BJ}{{\mathbb J}}

\def\be#1\ee{\begin{equation}#1\end{equation}}
\def\bea#1\eea{\begin{eqnarray}#1\end{eqnarray}}
\def\bean#1\eean{\begin{eqnarray*}#1\end{eqnarray*}}

\newcommand{\Tr}{\operatorname{\rm Tr}}

\newtheorem{definition}{Definition}[section]
\newtheorem{theorem}[definition]{Theorem}
\newtheorem{lemma}[definition]{Lemma}

\newtheorem{proposition}[definition]{Proposition}

\catcode `!=11

\newdimen\squaresize
\newdimen\thickness
\newdimen\Thickness
\newdimen\ll! \newdimen \uu! \newdimen\dd! \newdimen \rr! \newdimen
\temp!

\def\sq!#1#2#3#4#5{%
\ll!=#1 \uu!=#2 \dd!=#3 \rr!=#4
\setbox0=\hbox{%
 \temp!=\squaresize\advance\temp! by .5\uu!
 \rlap{\kern -.5\ll!
 \vbox{\hrule height \temp! width#1 depth .5\dd!}}%
 \temp!=\squaresize\advance\temp! by -.5\uu!
 \rlap{\raise\temp!
 \vbox{\hrule height #2 width \squaresize}}%
 \rlap{\raise -.5\dd!
 \vbox{\hrule height #3 width \squaresize}}%
 \temp!=\squaresize\advance\temp! by .5\uu!
 \rlap{\kern \squaresize \kern-.5\rr!
 \vbox{\hrule height \temp! width#4 depth .5\dd!}}%
 \rlap{\kern .5\squaresize\raise .5\squaresize
 \vbox to 0pt{\vss\hbox to 0pt{\hss $#5$\hss}\vss}}%
}
 \ht0=0pt \dp0=0pt \box0
}

\def\vsq!#1#2#3#4#5\endvsq!{\vbox to \squaresize{\hrule
width\squaresize height 0pt%
\vss\sq!{#1}{#2}{#3}{#4}{#5}}}

\newdimen \LL! \newdimen \UU! \newdimen \DD! \newdimen \RR!

\def\vvsq!{\futurelet\next\vvvsq!}
\def\vvvsq!{\relax
  \ifx     \next l\LL!=\Thickness \let\continue=\skipnexttoken!
  \else\ifx\next u\UU!=\Thickness \let\continue=\skipnexttoken!
  \else\ifx\next d\DD!=\Thickness \let\continue=\skipnexttoken!
  \else\ifx\next r\RR!=\Thickness \let\continue=\skipnexttoken!
  \else\def\continue{\vsq!\LL!\UU!\DD!\RR!}%
  \fi\fi\fi\fi
  \continue}

\def\skipnexttoken!#1{\vvsq!}

\def\place#1#2#3{\vbox to 0pt{\vss
\rlap{\kern#1\squaresize
  \raise#2\squaresize\hbox{$#3$}}
\vss}}

\catcode `!=12

\begin{document}
\maketitle


\tableofcontents


%
%





\vspace*{1cm}

The results in this paper form a step in the direction
of understanding the behavior of non-intersecting
Brownian motions on $\BR$ (Dyson's Brownian motions),
when the number of particles tends to $\iy$. Consider
$n$ Brownian particles leaving from points
$a_1<\ldots<a_p$ and forced to end up at $b_1<\ldots<
b_q$ at time $t=1$. It is clear that, when $n\rightarrow
\iy$, the equilibrium measure for $t\sim 0$ has its
support on $p$ intervals and for $t\sim 1$ on $q$
intervals. It is also clear that, when $t$ evolves,
intervals must merge, must disappear and be created,
leading to various {\em phase transitions}, depending on
the respective fraction of particles leaving from the
points $a_i$ and arriving at the points $b_j$. Therefore
the region ${\cal R}$ in the space-time strip $(x,t)$
formed by the support $(\subset \BR)$ of the equilibrium
measure as a function of time $0\leq t\leq 1$ will
typically present singularities of different types.

Near the moments, where a phase transition takes place,
one expects to find in the limit $n  { \nearrow} \iy$ an
infinite-dimensional diffusion, {\em a Markov cloud},
having some universality properties. Universality here
means that the infinite-dimensional diffusion is to
depend on the type of singularity only. These Markov
clouds are infinite-dimensional diffusions, which `in
principle' could be described by an infinite-dimensional
Laplacian with a drift term. We conjecture that each of
the Markov clouds obtained in this fashion is related to
some {\em integrable system}, which enables one to
derive a non-linear (finite-dimensional) PDE, satisfied
by the joint probabilities. The purpose of this paper is
to show, for a simple model leading to a cusp, that the
joint probabilities at different times, do satisfy such
a non-linear PDE. The interrelation between all such
equations and the ``initial" and ``final" ($t\rightarrow
\pm \iy$) conditions are interesting and challenging
open problems. Moreover, special cases have shown an
intimate connection between the integrable system and
the Riemann-Hilbert problem associated with the
singularity.

The {\bf first question} is the study of the finite
Brownian motion model, which, as will be explained in
section 2, hinges on the study of the coupled Gaussian
Hermitian random matrix ensemble ${\cal H}_n$ with
external source $A$, given by coupling terms
$c_1,\ldots,c_m$ and the diagonal matrix (set
$n=k_1+k_2$)
 \be
 A:= \left(\begin{array}{cccccc}
   \al\\
   &\ddots& & & &{\bf O}\\
   & & \al\\
   & & &-\al\\
   &{\bf O}& & &\ddots\\
   & & & & &-\al
\end{array}\right)\begin{array}{l}
\updownarrow k_1\\
\\
\\
\\
\updownarrow k_2
\end{array}.
  \label{0.1}\ee
    The probability
  of such an ensemble
  is defined by
\bean\lefteqn{\BP_n(\al,c_1,\ldots,c_m;E_1,\ldots, E_m)}\\
& :=&
 \frac{1}{Z_n}\int_{\prod^m_{\ell=1}{\cal
H}_n(E_{\ell})}\!\!\!\!\!\!\!e^{-\frac{1}{2}\Tr
(M_1^2+\ldots+M_m^2-2c_1M_1M_2-\ldots-2c_{m-1}M_{m-1}M_{m}-2AM_m)}
dM_1\ldots dM_m. \eean
  Given a disjoint union of intervals and the associated
  algebra of differential operators
  $$
  E_{\ell}:=
\bigcup^{r_\ell}_{i=1}[b^{(\ell)}_{2i-1},b^{(\ell)}_{2i}]\subset
\BR ~~\mbox{and}~~
 {\cal
   D}_k(E_{\ell})=\sum_{i=1}^{2r}(b_i^{(\ell)})^{k+1}\frac{\pl}{\pl
b_i^{(\ell)}},\qquad 1\leq\ell\leq m$$
 and given the tridiagonal matrix
\be
 J^{-1}:=
J^{-1}(c_1,\ldots,c_{m-1}):=\left(\begin{array}{ccccccccccc}
-1& &c_1\\
& & & \ddots& & {\bf O}\\
c_1& &-1& & \ddots& \\
 &\ddots& &\ddots & &\ddots\\
& & \ddots&  & -1& &c_{m-1}\\
 &{\bf O}& &\ddots\\
& & &  &  c_{m-1}& &-1
\end{array}
\right), \label{tridiag} \ee
 define the following differential operators:
  \bean
 \bar \AR_1^{\pm}
  &:=&-\frac{1}{2}\left(
  \sum_{j=1}^m J_{mj}   \DR_{-1}(E_{j})
   \pm\frac{\pl}{\pl \alpha}\right)
  \\
  \bar\CR_1&:=&\sum_{j=1}^m J_{1j} 
  \DR_{-1}(E_{j})
  \\
  \\
   \bar\AR_2
  &:=&\frac{1}{2}
 \left(
  \DR_0(E_{m})- \alpha\frac{\pl}{\pl \alpha}
   - c_{m-1}\frac{\pl}{\pl c_{m-1}}
  \right)
 \\
\bar{\cal C}_2&:=&-
  \DR_0(E_{1}) + c_1\frac{\pl}{\pl c_1}
 .\eean

\begin{theorem} 
%
%
%
  The log of the probability $\BP_n\left(\alpha;c_1,\ldots,c_{m-1};E_{1},\ldots,E_{m}\right)
$
   satisfies
  a  \underline {fourth-order} PDE in $\alpha$,
  $c_1,\ldots,c_{m-1}$
  and in the endpoints $b^{(\ell)}_1,..., b^{(\ell)}_{2r}$
  of the sets $E_{\ell}$, with
  \underline {quartic non-linearity}\footnote{in terms of the Wronskians
$\{f,g\}_X=gXf-fXg$.}:
 \bea
0&=& \Bigl(F^+\bar{\cal C}_{ 1}G^-+F^-\bar{\cal C}_{
1}G^+ \Bigr) \{F^-,F^+\}_{\bar{\cal C}_{ 1}}
 -
  \Bigl(F^+ G^- +F^- G^+
\Bigr) \bar{\cal C}_{ 1} \{F^-,F^+\}_{\bar{\cal C}_{ 1}} \no \\
&=& \det \left( \begin{array}{cccc}
  -G^+& \bar{\cal C}_1F^+  & - F^+ &0\\
   G^-& \bar{\cal C}_1F^-  & - F^- &0\\
   -\bar{\cal C}_1 G^+ &\bar{\cal C}_1^2F^+&0&-F^+\\
    \bar{\cal C}_1 G^- &\bar{\cal C}_1^2F^-&0&-F^-\\
  \end{array}\right)
\label{BasicEquation}\eea
where
\bean F^\pm
 &:=& \bar{\cal A}_1^\pm \bar{\cal C}_1
 \log \BP_{n}+k_{\left\{{1\atop 2}\right\}}J_{1m},
 \\ \\
 G^\pm&:=&\left\{
       \left(   \bar{\cal A}_2 \bar{\cal C}_1
 \pm  J_{1m}\frac{\pl}{\pl \alpha}\right)
  \log \BP_{n} \mp K_{\{{1\atop 2}\} }~,~
          F^\pm
       \right\}_{ \bar{\cal C}_{ 1}}
      \\ \\
      &&~~~~~~~~~~~~~~~~~~~~~+\left\{
         \left(  \bar{\cal
C}_2 \pm 2\alpha J_{1m} \bar{\cal C}_1\right){\cal
A}_1^\pm \log \BP_{n}~,~   F^\pm
     \right\}_{ \bar{\cal A}_{ 1}^\pm}
     \eean
     with $K_{\{{1\atop 2}\} }$ a constant,
     depending on $\al,c_i, k_1$ and $k_2$,
     $$
    K_{\{{1\atop 2}\} }:= J_{1m}\left(2k_{\{{1\atop 2}\}}\alpha
 J_{mm}-\frac{ k_1k_2}{\alpha} \right)
 .$$
%
%
\end{theorem}

Note the {\em robustness} of these equations: the
equations always have the same form
(\ref{BasicEquation}), regardless of the length of the
chain of matrices; only the quantities $F^{\pm}$,
$G^{\pm}$ and $H_i^{\pm}$ change, via some minors of the
matrix $J$.

 \bigbreak

The {\bf second question} concerns a simple model of
non-intersecting Brownian motions on $\BR$ and their
behavior, when the number of particles tends to $\iy$.

Consider $n=2k$ non-intersecting Brownian motions on
$\BR$,
   all starting at
the origin, such that the $k$ left paths end up at $-a$
and the $k$ right paths end up at $+a$ at time $t=1$;
see \cite{AptBleKui,TW,Okounkov,AvM4}.
 Inspired by \cite{Grabiner,Johansson,AptBleKui}, the
Karlin-McGregor formula \cite{Karlin} enables one to
express the transition probability $\BP^{\pm a}_0$ in
terms of the Gaussian Hermitian random matrices in a
chain $\BP_n(a;E)$ with external source, explained
above; this will be done in section 2.

Let now the number $n=2k$ of particles go to infinity,
and let the points $a$ and $-a$ go to $\pm \iy$. This
forces the left $k$ particles to $-\iy$ at $t=1$ and the
right $k$ particles to $+\iy$ at $t=1$. Since the
particles all leave from the origin at $t=0$, it is
natural to believe that for small times the equilibrium
measure (mean density of particles) is supported by {\em
one interval}, and for times close to $1$, the
equilibrium measure is supported by {\em two intervals}.
With a precise scaling, $t=1/2$ is critical in the sense
that for $t<1/2$, the equilibrium measure for the
particles is indeed supported by one, and for $t>1/2$,
by two intervals. The heart-shaped region $\cal R$
formed by the support of the equilibrium measure as a
function of time $0\leq t\leq 1$ has thus a cusp at
$t=1/2$.
 The {\em Pearcey process} ${\cal P}(s)$ is now
defined as the motion of an infinite number of
non-intersecting Brownian paths, just around time
$t=1/2$, with a precise scaling; see
\cite{AptBleKui,TW,Okounkov,AvM4}.  The joint
probability that the Pearcey process avoids the windows
$E_1,\ldots, E_m$ at times $s_1,\ldots,s_m$ is defined
by
%
\bean \lefteqn{ \BP\left(\begin{array}{c} \!\!{\cal
P}(s_1) \cap E_1= \emptyset\!\!
 \\ \vdots \\
 \!\!{\cal P}(s_m) \cap E_m= \emptyset\!\!
 \end{array}\right)  }\\
 &&
 := \lim_{z\rightarrow
0} \BP_0^{\pm 1/z^2}\left.\left(
 \begin{array}{c}
 \mbox{all}~x_j \bigl(\frac{1+s_{1}z^2}{2} \bigr) \notin
 zE_{1}\\
  \vdots
 \\
 \mbox{all}~x_j \bigl(\frac{1+s_{m}z^2}{2} \bigr) \notin
 zE_{m}
 \end{array}
 ; ~1\leq j\leq
n
  \right)\right|_{n=\frac{2}{z^4}}
  , \eean
where $\BP_0^{\pm a}$ was defined above. %
The main result of this paper is to show that the
infinite-dimensional diffusion equation for the Pearcey
process can be replaced by a finite-dimensional
non-linear PDE, which is intimately related to the
3-component KP hierarchy and which we now describe.


Given $E _{\ell}:=
\bigcup^{r_\ell}_{i=1}[x^{(\ell)}_{2i-1},x^{(\ell)}_{2i}]\subset
\BR$, define the \underline{space} and \underline{time}
gradients
%
%
$$
   {\cal
X}_{-1}:=\sum_{\ell=1}^m \sum_{i=1}^{2r_\ell}
\frac{\pl}{\pl x_i^{(\ell)}},~~~{\cal
T}_{-1}=\sum_{\ell=1}^m  \frac{\pl}{\pl s_{\ell}},
 $$
 \underline{space} and \underline{time} {\em Euler operators}
 $\tilde{\cal
X}_{0}$ and $\tilde{\cal T}_{0}$
 and a mixed \underline{space-time} operator
 $\tilde{\cal
X}_{-1}$,
$$
   {\cal
X}_{0}:=\sum_{\ell=1}^m \sum_{i=1}^{2r_\ell}
 x_i^{(\ell)} \frac{\pl}{\pl x_i^{(\ell)}},~~~{\cal
T}_{0}=\sum_{\ell=1}^m  s_{\ell}\frac{\pl}{\pl
s_{\ell}},~~~ \tilde{\cal X}_{-1}=\sum_{\ell=1}^m
 s_{\ell}\sum_{i=1}^{2r_\ell}\frac{\pl}{\pl
 x_i^{(\ell)}}.
$$

\bigbreak

\begin{theorem}
Then $$\BQ(s_1,\ldots,s_m;E_{1},\ldots,E_{m}):=\log
\BP\left(\begin{array}{c} \!\!{\cal P}(s_1) \cap E_1=
\emptyset\!\!
 \\ \vdots \\
 \!\!{\cal P}(s_m) \cap E_m= \emptyset\!\!
 \end{array}\right)  $$
   satisfies a 4th order
and 3rd degree PDE, which can be written as a single
Wronskian in the gradient ${\cal X}_{-1}$:
\bea
\lefteqn{\Biggl\{{\cal X}^2_{-1}{\cal T}_{-1}\BQ~,~
  \frac{1}{8}\left\{ {\cal X}_{-1}{\cal T}_{-1}\BQ ~
,~{\cal X}^2_{-1} \BQ\right\}_{{\cal X}_{-1}}}
\label{PearceyEquation} \\
&& ~~~ + \left({\cal X}_0\!+\!2{\cal T}_0\!-\!2 \right)
{\cal X}_{-1}^2\BQ - 4(\tilde{\cal X}_{-1}{\cal
X}_{-1}-{\cal T}^2_{-1}) {\cal T}_{-1}\BQ
\Biggr\}_{{\cal X}_{-1}}=0. \no \eea
In particular $
 \BQ(s;E)=\log\BP\left({\cal P}(s ) \cap E =
\emptyset\right) $ satisfies
  $$
 \left\{
  {\cal X}_{-1}^2\frac{\pl \BQ}{\pl
s} ~,~
  \frac{1}{8}\Bigl\{
{\cal X}_{-1}\frac{\pl \BQ}{\pl s}~,~{\cal
X}_{-1}^2\BQ\Bigr\}_{{\cal X}_{-1}} +  ({\cal
X}_0-2){\cal X}_{-1}^2\BQ
  +4
  \frac{\pl^3 \BQ}{\pl s^3}  \right\}_{{\cal X}_{-1}}=0.
$$

\end{theorem}
Notice here as well the robustness of the equations. The
shape of the equation (\ref{PearceyEquation}) is the
same, regardless of the number of times one considers.
Moreover the equations are ``{\em commutative}": the
times and windows can be permuted simultaneously. Notice
that the term containing $\tilde {\cal X}_{-1}$ is the
only one which ties up the time $s_{i}$ with the precise
set $E_i$. We expect that this equation can be used to
derive large-time asymptotics, when $t\rightarrow \pm
\iy$. Also one expects that the PDE's for the sine
 and Airy processes \cite{AvM3} can be obtained from this
equation by an appropriate scaling limit. These
questions remain challenging open problems.

\newpage

\section{Gaussian Hermitian random matrices coupled in a
chain with external source}

The present paper studies $m$ Gaussian Hermitian random
matrices $M_i \in {\cal H}_n$, coupled in a chain with
external source $A$, given by the diagonal matrix (set
$n=k_1+k_2$)
 \be
 A:= \left(\begin{array}{cccccc}
   \alpha\\
   &\ddots& & & &{\bf O}\\
   & &\alpha\\
   & & &-\alpha\\
   &{\bf O}& & &\ddots\\
   & & & & &-\alpha
\end{array}\right)\begin{array}{l}
\updownarrow k_1\\
\\
\\
\\
\updownarrow k_2
\end{array},
  \label{0.1}\ee
  and given by the coupling terms $c_1,\ldots,c_{m-1}$;
  its density is given by
  \be
  \frac{1}{Z_n} 
  e^{-\frac{1}{2}\Tr
(M_1^2+\ldots+M_m^2-2c_1M_1M_2-\ldots-2c_{m-1}M_{m-1}M_{m}-2AM_m)}
dM_1\ldots dM_m
  \ee
    For each index $1\leq\ell\leq m$, consider
      a disjoint union of intervals
  $
  E_{\ell}:=
\bigcup^r_{i=1}[b^{(\ell)}_{2i-1},b^{(\ell)}_{2i}]\subset
\BR$, and define the associated algebra of differential
operators
  \be{\cal
   D}_k{(E_\ell)}=\sum_{i=1}^{2r}(b_i^{(\ell)})^{k+1}\frac{\pl}{\pl
b_i^{(\ell)}},\qquad 1\leq\ell\leq m.\ee
Consider the following probability:

\bea
\lefteqn{\BP_n\left(\alpha;c_1,\ldots,c_{m-1};E_{1},\ldots,E_{m}\right)
}\no\\ \no\\
 &:=&
   \BP\left((M_1, \ldots, M_m)\in\prod^m_{\ell=1}{\cal
H}_n~,~\mbox{with}\left\{\begin{array}{c}\mbox{\!$M_1$-spectrum
in
$E_1$}\\
 \vdots\\
 \mbox{$M_{m}$-spectrum in
$E_m$ }\end{array}\right\}
\right) \no\\
  &=&\!\!
  \frac{1}{Z_n'}\int_{\prod^m_{\ell=1}{\cal
H}_n(E_{\ell})}\!\!\!\!\!\!\!e^{-\frac{1}{2}\Tr
(M_1^2+\ldots+M_m^2-2c_1M_1M_2-\ldots-2c_{m-1}M_{m-1}M_{m}-2AM_m)}
dM_1\ldots dM_m \no\\
 \eea
To be clear, this integral is to be taken over the space
of $m$-uples of Hermitian matrices, with
\mbox{$M_1$-spectrum in $E_1,\ldots,~M_{m}$-spectrum in
$E_m$ }, and $Z_n'$ is the above integral with all the
$E_i$ replaced by $\BR$.


\begin{proposition}.
The following holds\footnote{Throughout this paper,
$\Dt_n(x)=\prod_{1\leq i< j\leq n}(x_i-x_j)$ is the
Vandermonde.}:
\bean
\lefteqn{\BP_n\left(\alpha;c_1,\ldots,c_{m-1};E_1,\ldots,E_m\right)
}\\
& =& \frac{1}{Z_n}\int_{\prod_1^m E_i^{k_1+k_2}}
\Dt_{k_1+k_2}(y^{(1)})
 \\
   && \Dt_{k_1}(y^{(m)'})
\prod_{i=1}^{k_1}e^{-\frac{1}{2}\sum_{\ell=1}^{m}
y_i^{(\ell)2}
  + \sum_{\ell=1}^{m-1}c_{\ell}y_i^{(\ell)}y_i^{(\ell+1)}
   +\alpha y_i^{(m)}  }
     \prod_{\ell=1}^{m}dy_i^{(\ell)}
      \\
 &&   \Dt_{k_2}(y^{(m)''})
  \!\!\prod_{i=k_1+1}^{k_1+k_2}\!\!
   e^{-\frac{1}{2}\sum_{\ell=1}^{m}
y_i^{(\ell)2}
  + \sum_{\ell=1}^{m-1}c_{\ell}y_i^{(\ell)}y_i^{(\ell+1)}
   -\alpha y_i^{(m)}}
    \prod_{\ell=1}^{m}dy_i^{(\ell)}
    ,\eean
    where
    $$
y^{(m)'}:=(y_{1}^{(m)},\ldots,y_{k_1}^{(m)})
~~~~\mbox{and}~~~~
y^{(m)''}:=(y_{k_1+1}^{(m)},\ldots,y_{k_1+k_2}^{(m)})
.$$

   \end{proposition}

The {\em proof} of this statement is a standard
application of the Harish-Chandra-Bessis-Itzykson-Zuber
formula, combined with the techniques explained in the
next section.


\section{Non-intersecting Brownian motions}

Consider $n=k_1+k_2$ non-intersecting Brownian motions
on $\BR$ (Dyson's Brownian motions),
   all starting at
the origin, such that the $k_2$ left paths end up at
$-a$ and the $k_1$ right paths end up at $+a$ at time
$t=1$:
 \bea
\lefteqn{ \BP_0^{\pm a}\left(\mbox{all $x_j(t)\in
E$}\right)
  }
   \no\\
&:=& \BP \left(\begin{tabular}{c|c}
   & all $x_j(0) =0$\\
  all $x_j(t)\in E$  &
   $k_2$ left paths end up at $-a$ at time $t=1$,\\
   &
$k_1$ right paths end up at $+a$ at time $t=1$
\end{tabular}\right)\no\\
\eea
%
%
In the Proposition below we shall be using the
Karlin-McGregor formula for non-intersecting Brownian
motions $x_j(t)$ for $0<t<1$:

\bea
 \lefteqn{\BP\left( \mbox{all}~ x_i(t)
  \in E,~1\leq i\leq n \left|
  \begin{array}{l}
  \mbox{given $x_i(0)=\gamma_i$}\\
  \mbox{given $x_i(1)=\delta_i$}
   \end{array}
  \right.\right)}\no\\
  &=&
  \int_{E^n}\frac{1}{Z_{n}}\det(p(t;\gamma_i,x_j))_{1\leq i,j\leq
n}\det(p(1-t;x_{i'},\delta_{j'}))_{1\leq i',j'\leq n}
 \prod_1^n dx_i\no\eea
  for
  \be
  p(t,x,y):= \frac{1}{\sqrt{ \pi t}
  }~e^{-\frac{(y-x)^2}{ t}}.
  \label{BrownianTransition}\ee
Consider now the Brownian motions at different times
$$
0=t_0<t_1<t_2<\ldots <t_{m-1}<t_m<t_{m+1}=1
$$
 and set
 $$
\tau_i=t_{i+1}-t_i
\mbox{~~and~~}
\frac{1}{\sg_j}=\frac{1}{t_j-t_{j-1}}+\frac{1}{t_{j+1}-t_j},\quad\mbox{for~}0\leq
j\leq m.
$$
Considering
      $m$ disjoint unions of intervals
  $
  E_{\ell}:=
\bigcup^r_{i=1}[\tilde b^{(\ell)}_{2i-1},\tilde
b^{(\ell)}_{2i}]\subset \BR$ for $1\leq \ell \leq m$, we
show that the two probabilities $ \BP_0^{\pm a}$ and
$\BP_n$, as in (1.4) and (2.1), are related by a mere
change of variables:

\begin{proposition}

For $0=t_0<t_1<t_2<\ldots <t_{m-1}<t_m<t_{m+1}=1$,
\bean
  {\BP_0^{\pm a}(\mbox{all~}x_i(t_1)\in\tilde
E_1,\ldots,\mbox{all~}x_i(t_m)\in\tilde E_m)}
&=&\BP_n\left(\alpha;c_1,\ldots,c_{m-1};E_1,\ldots,E_m\right)
\eean
 upon setting
\be
E_{\ell}
=
 \tilde E_\ell\sqrt{\frac{2(t_{\ell+1}-t_{\ell-1})}
  {(t_{\ell+1}-t_{\ell})(t_{\ell}-t_{\ell-1})}}
 ~
,~c_j
 =\sqrt{ \frac{(t_{j+2}-t_{j+1})(t_{j}-t_{j-1})}
   {(t_{j+2}-t_{j})(t_{j+1}-t_{j-1})}}\label{changeof
   variables1}\ee
    \mbox{and}
   \be\al
=a\sqrt{\frac{2(t_{m}-t_{m-1})}{(1-t_m)(1-t_{m-1})}}
  .
 \label{changeof
   variables2} \ee 
%
%

\end{proposition}

\proof In the following computation, we shall be using
the notation $$
x^{(m)'}=(x_1^{(m)},\ldots,x_{k_1}^{(m)}),\quad
x^{(m)''}=(x^{(m)}_{k_1+1},\ldots,x^{(m)}_{k_1+k_2})
 .$$
 Remembering $p(t,x,y)$ is the Brownian transition
 probability (\ref{BrownianTransition}),
one computes:
\bean
\lefteqn{\BP_0^{\pm a}(\mbox{all~}x_i(t_1)\in\tilde
E_1,\ldots,\mbox{all~}x_i(t_m)\in\tilde
E_m)}\\
&=&\lim_{\begin{array}{c}
{\scriptstyle x_i^{(0)}\rg 0}\\
{\scriptstyle a_1,\ldots,a_{k_1}\rg  a}\\
{\scriptstyle a_{k_1+1},\ldots,a_n\rg -a}
\end{array}}\frac{1}{Z_n}\int_{\tilde E^n_1\times\ldots\times\tilde
E_m^n}\prod_{1\leq i\leq  n\atop{1\leq j\leq
m}}dx_i^{(j)}\\
\\
& & \det\left(p(t_1,x_i^{(0)},x_j^{(1)})\right)_{1\leq
i,j\leq
n}\det\left(p(t_2-t_1,x_i^{(1)},x_j^{(2)})\right)_{1\leq
i,j\leq n}\\
\\
& &\ldots\det \left(p(t_m-t_{m-1},x_i^{(m-1)},x_j^{(m)})
\right)_{1\leq i,j\leq n}
\det\left(p(1-t_m,x_i^{(m)},a_j)\right)_{1\leq
i,j\leq n}\\
\\
&=&\lim_{\begin{array}{c}
{\scriptstyle a_1,\ldots,a_{k_1}\rg  a}\\
{\scriptstyle a_{k_1+1},\ldots,a_n\rg -a}
\end{array}}\frac{1}{Z'_n}
 \int_{\tilde E^n_1\times\ldots\times\tilde E_m^n}
\prod_{1\leq i\leq  n\atop{1\leq j\leq
m}}dx_i^{(j)}\Dt_n(x^{(1)})\prod^m_{j=1}
e^{-\sum^n_{i=1}x_i^{(j)2}\left(\frac{1}{t_j-t_{j-1}}
+\frac{1}{t_{j+1}-t_{j}}\right)}
\\
\\
& &\det\left(
e^{\frac{2x_i^{(1)}x_j^{(2)}}{t_2-t_1}}\right)_{1\leq
i,j\leq n}\ldots\det\left(
e^{\frac{2x_i^{(m-1)}x_j^{(m)}}{t_m-t_{m-1}}}\right)_{1\leq
i,j\leq n}\det\left(e^{\frac{2x_i^{(m)}a_j}{1-t_m}}\right)_{1\leq i,j\leq n}\\
\\
&\stackrel{*}{=}&\frac{(n!)^m}{Z''_n}\int_{\tilde
E^n_1\times\ldots\times\tilde E_m^n}\prod_{1\leq i\leq
n\atop{1\leq j\leq m}}dx_i^{(j)}\Dt_n(x^{(1)})
 \\
 &&\prod_{i=1}^{k_1}(x_i^{(m)})^{i-1}e^{-\sum^m_{j=1}\frac{x_
i^{(j)2}}{\sg_j}
+\sum_{j=1}^{m-1}\frac{2x_i^{(j)}x_i^{(j+1)}}{t_{j+1}-t_j}
 +\frac{2ax_i^{(m)}}{1-t_m}}\
\
\\
&
&\prod^{k_1+k_2}_{i=k_1+1}(x_i^{(m)})^{i-k_1-1}e^{-\sum^m_{j=1}
\frac{x_i^{(j)
2}}{\sg_j}+\sum^{m-1}_{j=1}\frac{2x_i^{(j)}x_i^{(j+1)}}{t_{j+1}
-t_j}
-\frac{2ax_i^{(m)}}{1-t_m}}\\
\\
&\stackrel{**}{=}&\frac{(n!)^m}{(k_1!)(k_2)!Z''_n}\int_{\tilde
E^n_1\times\ldots\times\tilde E_m^n}\Dt_n(x^{(1)})\\
\\
& &\Dt_{k_1}(x^{(m)'})\prod^{k_1}_{i=1}
 e^{-\sum^m_{j=1}\frac{x_i^{(j)2}}{\sg_j}+\sum_{j=1
}^{m-1}\frac{2x_i^{(j)}x_i^{(j+1)}} {t_{j+1}-t_j}
+\frac{~2ax_i^{(m)}}{1-t_m}
}\prod^m_{j=1}dx_i^{(j)}\\
\\
&
&\Dt_{k_2}(x^{(m)''})\prod^{k_1+k_2}_{i=k_1+1}e^{-\sum^m_{j=1}
\frac{x_i^{(j)
2}}{\sg_j}+\sum_{j=1}^{m-1}\frac{2x_i^{(j)}x_i^{(j+1)}}
{t_{j+1}-t_j} -\frac{~2ax_i^{(m)}}{1-t_m}
}\prod^m_{j=1}dx_i^{(j)} \eean
\bean &=&\frac{(n!)^m}{Z''_n k_1!k_2!}\int_{\prod_1^m
E_i^n}\prod_{1\leq i\leq n\atop{1\leq j\leq
m}}dy_i^{(j)}\Dt_n(y^{(1)})
\\
\\
& &\Dt_{k_1}(y^{(m)'})\prod^{k_1}_{i=1}
e^{-\frac{1}{2}\sum^m_{j=1}y_i^{(j)2}
   +\sum^{m-1}_{j=1}c_jy_i^{(j)}y
_i^{(j+1)}
+\al y_i^{(m)}}\\
\\
&
&\Dt_{k_2}(y^{(m)''})\prod^{k_1+k_2}_{i=k_1+1}e^{-\frac{1}{2}\sum^m_{j=1}y_i^{(j)2}+\sum^{m-1}_{j=1}c_jy_
i^{(j)}y_i^{(j+1)} -\al y_i^{(m)}} \eean
by setting, for $1\leq i \leq n,~ 1\leq j \leq m-1$,
$$
\frac{x_i^{(\ell)}}{\sqrt{\sg_{\ell}}}=\frac{y_i^{(\ell)}}{\sqrt{2}},\qquad
c_j=\frac{\sqrt{\sg_j\sg_{j+1}}}{t_{j+1}-t_j}\mbox{~and~}\al
=\frac{\sqrt{2\sg_m }}{1-t_m}a ;
$$
 the change of variables $x^{(\ell)}\mapsto y^{(\ell)}$
 induces a change of variables for the boundary terms of
 the integrals
$$E_{ \ell }=\tilde
E_{\ell}\sqrt{\frac{2}{\sg_{\ell}}},$$ thus confirming
(\ref{changeof
   variables1}) and (\ref{changeof
   variables2}).



Identity $\stackrel{*}{=}$ in the previous set of
identities is established by means of the following
argument, which we explain for the indices $j=1,2$:

\bean \lefteqn{\prod_{1\leq i\leq n\atop{1\leq j\leq
2}}dx_i^{(j)}\Dt _n(x^{(1)})\det\left(
e^{\frac{2x_i^{(1)}x_j^{(2)}}{\tau_1}}\right)_{1\leq
i,j\leq
n}\det\left(e^{\frac{2x_i^{(2)}x_j^{(3)}}{\tau_2}}\right)_{1\leq
i,j\leq n}
}\\
\\
& =&\prod_{1\leq i\leq n\atop{1\leq j\leq
2}}dx_i^{(j)}\Dt
_n(x^{(1)})\sum_{\pi}(-1)^{\pi}\prod^n_1
e^{\frac{2x^{(1)}_{\pi (i)}x_i^{(2)}}{\tau_1}}
 \sum_{\pi'}(-1)^{\pi'}\prod^n_1
e^{\frac{2x^{(2)}_{\pi'(j)}x_j^{(3)}}{\tau_2}}\\
\eean
Upon setting
\bean
x^{(2)}_{\pi'(j)}&\longmapsto&x_j^{(2)},\\
\\
x^{(1)}_{\pi(\pi'(i))}&\longmapsto&x_i^{(1)}, \eean
this expression turns into
\bean\\
& {\longmapsto}&\prod_{1\leq i\leq n\atop{1\leq j\leq
2}}dx_i^{(j)}\sum_{\pi,\pi'}(-1)^{\pi+\pi'}\Dt
_n(x^{(1)}_{(\pi\circ\pi')^{-1}(1)},\ldots,
 x^{(1)}_{(\pi\circ\pi')^{-1}(n) })\prod_1^n
e^{\frac{2x_i^{(1)}x_i^{(2)}}{\tau_1}}
 \prod^n_1e^{\frac{2x_j^{(2)}x_j^{(3)}}{\tau_2}}\\
\\
&=&(n!)^2\Dt_n(x^{(1)})\prod_{1\leq i\leq n\atop{1\leq
j\leq 2}}dx_i^{(j)}
\prod_1^ne^{\frac{2x_i^{(1)}x_i^{(2)}}{\tau_1}+\frac{2x_i^{(2)}x_i^{(3)}
}{\tau_2}}. \eean
Then one uses the symmetry of the integration ranges
 vis-\`a-vis $i$. In general, one makes not $2$ but $m$
 synchronized changes of variables.

 Identity $\stackrel{**}{=}$ follows by simultaneously
 setting
 \bean
  x_{i}^{(\ell)}&\mapsto& x_{\pi'(i)}^{(\ell)},
  ~~1\leq i\leq k_1,~1\leq \ell\leq m, ~\pi'\in S_{k_1} \\
  x_{k_1+j}^{(\ell)}&\mapsto& x_{k_1+\pi''(j)}^{(\ell)},
  ~~1\leq j\leq k_2,~1\leq \ell\leq m, ~\pi''\in S_{k_2}
  ,\eean
   subsequently summing over $\pi'\in S_{k_1},~\pi''\in
  S_{k_2}$, giving rise to $\Dt_{n_1}(x^{(m)'})$ and
  $\Dt_{n_2}(x^{(m)''}  )$, because of the presence of $\Dt_n(x^{(1)})$, and then dividing by $k_1!k_2!$,
 thus ending the proof of Proposition 2.1. \qed


When taking the limit for $n\rightarrow \iy$, we shall
need the following scaling (see \cite{AptBleKui}),
assuming $k=k_1=k_2$:
 \be
 n=2k=\frac{2}{z^4},\quad \pm a=\pm\frac{ 1}{z^2},\quad
 x_i\mapsto x_iz,\quad
 t_j=\frac{1}{2}(1+s_jz^2),~~~\mbox{~for~}z\rg 0.
 \label{0.10}
 \ee


  \begin{proposition}. Given $\tilde E_{\ell}=\bigcup^r_{i=1}[\tilde
b^{(\ell)}_{2i-1},\tilde b^{(\ell)}_{2i}]\subset\BR$,
the following holds:

$$
\BP_0^{\pm a}(\mbox{all~}x_i(t_1)\in\tilde
E_1,\ldots,\mbox{all~}x_i(t_m)\in\tilde
E_m)\Bigg\vert_{\begin{array}{l}
t_j=\frac{1}{2}(1+s_jz^2)\\
\tilde b_j^{(\ell)}=u_j^{(\ell)}z \\
a=\frac{1}{z^2}\\
n=\frac{2}{z^4}
\end{array}}
$$

   \be= \BP_n\left(\alpha;c_1,\ldots,c_{m-1};b^{(1)},\ldots,b^{(m)}\right)
  \Bigr|_{n=\frac{2}{z^4}} ,
   \label{scaling}\ee
   with
\bean
c_j=\frac{\sqrt{\sg_j\sg_{j+1}}}{t_{j+1}-t_j}&=&\sqrt{\frac{(t_j-t_{j-1})(t_{j+2}-t_{
j+1})}{(t_{j+1}-t
_{j-1})(t_{j+2}-t_j)}}\\
\\
&=&\left\{\begin{array}{ll}
\sqrt{\frac{(1+s_1z^2)(s_3-s_2)}{(1+s_2z^2)(s_3-s_1)}}&\mbox{for~}j=1\\
\\
\sqrt{\frac{(s_j-s_{j-1})(s_{j+2}-s_{j+1})}{(s_{j+1}-s_{j-1})
(s_{j+2}-s_j)}}&\mbox{for~}2\leq j\leq
m-2\\
\\
\sqrt{\frac{(s_{m-1}-s_{m-2})(1-s_mz^2)}{(s_m-s_{m-2})(1-s_{m-1}z^2)}}&\mbox{for~}j=m
-1
\end{array}\right.\\
\\
\al=\frac{a\sqrt{2\sg_m}}{1-t_m}&=&\frac{ \sqrt{2}}{z^2}
\frac{\sqrt{t_m-t_{m-1}}}{\sqrt{(1-t_m)(1-t
_{m-1})}}=\frac{2}{z}\frac{\sqrt{s_m-s_{m-1}}}{\sqrt{(1-s_mz^2)(1-s
_{m-1}z^2)}} \eean

\bean
 b_i^{(\ell)}=\tilde
b_i^{(\ell)}\sqrt{\frac{2}{\sg_{\ell}}}
 &=&\tilde
b_i^{(\ell)}
\sqrt{\frac{2(t_{\ell+1}-t_{\ell-1})}{(t_{\ell}-t_{\ell-1})
(t_{\ell+1}-t_{\ell})}} \\
&=& \left\{\begin{array}{l}
2u_i^{(1)}\sqrt{\frac{1+s_2z^2}{(1+s_1z^2)(s_2-s_1)}}
 ,~\mbox{~~for $\ell=1$}\\
 2u_i^{(\ell)}
 \sqrt{\frac{s_{\ell+1}-s_{\ell-1}}{(s_{\ell}-s_{\ell-1})
 (s_{\ell+1}-s_{\ell})}}
  ,~\mbox{~~for $2\leq \ell\leq m-1$} \\
2u_i^{(m)}\sqrt{\frac{1-s_{m-1}z^2}{(s_m-s_{m-1})(1-s_mz^2)}}
 ,~\mbox{~~for $\ell=m$}
\end{array}\right.
\eean

\end{proposition}

\proof Straightforward from Proposition 2.1, combined
with the scaling, appearing in (\ref{scaling}).\qed

 \newpage

\section{The inverse of a tridiagonal matrix and its derivatives}





Consider the $(k+1)\times (k+1)$ tridiagonal matrix,
with non-diagonal entries $c_1,\ldots,c_k$:

$$
J^{-1}(c_1,\ldots,c_k)=\left(\begin{array}{ccccccccccc}
-1& &c_1\\
& & & \ddots& & {\bf O}\\
c_1& &-1& & \ddots& \\
 &\ddots& &\ddots & &\ddots\\
& & \ddots&  & -1& &c_k\\
 &{\bf O}& &\ddots\\
& & &  &  c_k& &-1
\end{array}
\right)
$$
with
\be D(c_1,\ldots,c_k):=\left\{\begin{array}{ll}
\det J^{-1}(c_1,\ldots,c_k),&\mbox{for~}k\geq 1\\
\\
-1&\mbox{for~}k=0\\
\\
1&\mbox{for~}k=-1.
\end{array}\right.
\label{A} \ee Then one checks, that for $1\leq j\leq m$,

$$
J_{1j}(c_1,\ldots,c_{m-1})=(-1)^{j-1}c_1\ldots
c_{j-1}\frac{D(c_{j+1},\ldots,c_{m-1})}{D(c_1,\ldots,c_{m-1})}
$$
and
\be J_{mj}(c_1,\ldots,c_{m-1})=(-1)^{m-j}c_j\ldots
c_{m-1}\frac{D(c_1,\ldots,c_{j-2})}{D(c_1,\ldots,c_{m-1})}.
\label{B} \ee Define
\be
\begin{array}{lll}
r_j&:=\frac{J_{mj}}{J_{1j}}&\mbox{for~}1\leq j\leq m \\
 \\
&=0&\mbox{for~}j=0.
\end{array}\label{C}
\ee
For later use, we shall need the following identities:

\begin{lemma}

 \bean c_1\frac{\pl}{\pl c_1}\frac{J_{m1}J_{mi}}{J_{1i}}
  =-2J_{m1}^2 &\mbox{,}&  c_1\frac{\pl}{\pl
  c_1}
  \log \frac{J_{1m} }{J_{1i}}
  =\dt_{1i} ~,~c_1 \frac{\pl}{\pl c_{1}} {J_{mm}}
  =
  -2J_{m1}^2
   \\
  \\
  c_{m-1}\frac{\pl}{\pl c_{m-1}}\log\frac{J_{m1}J_{mi}}{J_{1i}}
  &=& 2\left(-1-2J_{mm}+\frac{J_{m1}J_{mi}}{J_{1i}}\right)
  \\
    c_{m-1}\frac{\pl}{\pl c_{m-1}}\log\frac{J_{m1} }{J_{1i}}
  &=&-1-\frac{2D(c_{i+1},\ldots,c_{m-2})}
   {D(c_{i+1},\ldots,c_{m-1})}=-1-
2\frac{J_{mm}J_{1i}-J_{m1}J_{mi}}{J_{1i}}
   \\
    c_{m-1}\frac{\pl}{\pl c_{m-1}}\log J_{mi}
  &=&
    1-2c_{m-1}J_{m,m-1}-\dt_{im}=-2J_{mm}-1-\dt_{im}
   \eean

  \end{lemma}

\proof At first notice that
\bea
 c_1\frac{\pl}{\pl
 c_1}D(c_1,\ldots,c_j)&=&-2c_1^2D(c_3,\ldots,c_j)\no\\
  &=&
 2\left(D(c_1,\ldots,c_{j})+D(c_2,\ldots,c_{j})\right)
\no \\ \no\\ \no\\
c_{m-1}\frac{\pl}{\pl c_{m-1}}D(c_j,\ldots,c_{m-1})&=&
 -2c_{m-1}^2D(c_j,\ldots,c_{m-3})\no\\
 &=&
 2\bigl(D(c_j,\ldots,c_{{m-1}})+D(c_j,\ldots,c_{m-2})\bigr)
 \no\\ \label{3.P},\eea
and {\footnotesize $$
 D(c_2,\ldots,c_{m-1})D(c_1,\ldots,c_{i-2})
 -D(c_1,\ldots,c_{m-1})D(c_2,\ldots,c_{i-2})
  = D(c_{i+1},\ldots,c_{m-1})\prod _1^{i-1}c_k^2
.$$}

 {\footnotesize \be
 D(c_1,\ldots,c_{m-2})D(c_{i+1},\ldots,c_{m-1})
 -D(c_1,\ldots,c_{m-1})D(c_{i+1},\ldots,c_{m-2})
  = D(c_{ 1},\ldots,c_{i-2})\prod _i^{m-1}c_k^2
.\label{3.Q}\ee}
Then one checks, by (\ref{B}) and (\ref{3.Q}), that
\bea \frac{J_{m1}J_{mi}}{J_{1i}} &=&\frac{c_i^2\ldots
c_{m-1}^2D(c_{1},\ldots,c_{i-2})}{D(c_{i+1},\ldots,c_{m-1})
 D(c_{1},\ldots,c_{m-1})}
\no\\
&=&\frac{D(c_{1},\ldots,c_{m-2})}{D(c_{1},\ldots,c_{m-1})}
-
\frac{D(c_{i+1},\ldots,c_{m-2})}{D(c_{i+1},\ldots,c_{m-1})}
 \no\\& =&J_{mm}-
  \frac{D(c_{i+1},\ldots,c_{m-2})}
  {D(c_{i+1},\ldots,c_{m-1})}.
  \label{C1}\eea
 %
 Hence
$$
\frac{D(c_{i+1},\ldots,c_{m-2})}{D(c_{i+1},\ldots,c_{m-1})}
=\frac{J_{mm}J_{1i}-J_{m1}J_{mi}}{J_{1i}},
$$
 and by (\ref{B}),
  \be
\frac{J_{m1} }{J_{1i}} =\frac{(-1)^{m-i}c_i \ldots
c_{m-1}
 }{D(c_{i+1},\ldots,c_{m-1})}
.\label{C2}\ee
 Moreover, explicit differentiation of (\ref{C1}) and using
 the identities (\ref{3.Q}), (\ref{B}), one is led to
 \bean
 \lefteqn{ c_1\frac{\pl}{\pl c_1}\frac{J_{m1}J_{mi}}{J_{1i}}
   }\\
   &=&
    \frac{-2c_i^2\ldots c_{m-1}^2
     \left(D(c_2,\ldots,c_{m-1})D(c_1,\ldots ,c_{i-2})-
 D(c_2,\ldots,c_{i-2})D(c_1,\ldots,c_{m-1})
 \right)}
     {D(c_{i+1},\ldots,c_{m-1})D^2(c_{1},\ldots,c_{m-1})}
\\
&=&
\frac{-2\prod_1^{m-1}c^2_i}{D^2(c_{1},\ldots,c_{m-1})}
\\
 &=&
 -2J_{m1}^2
 ,\eean
 and, setting $i=m$, yields at once the last identity on
 the first line of the statement of Lemma 3.1. Also by
 (\ref{C1}) and (\ref{3.Q})
{\footnotesize \bean
 \lefteqn{\frac{1}{2} c_{m-1}\frac{\pl}{\pl c_{m-1}}\log\frac{J_{m1}J_{mi}}{J_{1i}}
   }\\
   &=&\!\!
  -1-2J_{mm}+\frac{D(c_1,\ldots,c_{m-2})
  D(c_{i+1},\ldots,c_{m-1})-
   D(c_1,\ldots,c_{m-1})D(c_{i+1},\ldots,c_{m-2})}
   {D(c_1,\ldots,c_{m-1})D(c_{i+1},\ldots,c_{m-1})}
    \\
   &=& \!\!   -1-2J_{mm}+\frac{J_{m1}J_{mi}}{J_{1i}}
   \eean
}
From (\ref{C2}) compute at once
$$
c_1\frac{\pl}{\pl c_1}\log \frac{J_{1m}}{J_{1i}}
 =
 \dt_{1i}
 $$
 and from (\ref{C2}), (\ref{3.P}) and (\ref{C1}),
\bean
  c_{m-1}\frac{\pl}{\pl c_{m-1}}\log
\frac{J_{1m}}{J_{1i}}
 &=&
 1-2\frac{D(c_{i+1},\ldots,c_{m-1})+D(c_{i+1},\ldots,c_{m-2})}
  {D(c_{i+1},\ldots,c_{m-1})}
  \\
  &=&-1-2\frac{D(c_{i+1},\ldots,c_{m-2})}{D(c_{i+1},\ldots,c_{m-1})}
 \\
 &=&-1-
2\frac{J_{mm}J_{1i}-J_{m1}J_{mi}}{J_{1i}}.
 \eean
%
%
%
 Also, using
$$
D(c_1,\ldots,c_{m-1})=-D(c_1,\ldots,c_{m-2})-c_{m-1}^2D(c_1,\ldots,c_{m-3})
$$
and the explicit formula (\ref{B}) for $J_{mi}$, one
computes, using (\ref{3.P}),
\bean
c_{m-1}\frac{\pl}{\pl c_{m-1}}\log J_{mi}
&=&c_{m-1}\frac{\pl}{\pl
c_{m-1}}\log\left((-1)^{m-i}c_i\ldots
c_{m-1}\frac{D(c_1,\ldots,c_{i-2})}{D(c_1,\ldots,c_{m-1})}\right)\\
\\
&=& 1+2c_{m-1}^2\frac{D(c_1,\ldots,c_{m-3})}
{D(c_1,\ldots,c_{m-1})}  -\dt_{im}
\\
&=& 1 -2c_{m-1}J_{m,m-1} -\dt_{im} \\
&=&
 2\left(1-c_{m-1}J_{m,m-1}\right)-1 -\dt_{im}\\
&=& 2\left(\frac{D(c_1,\ldots,c_{m-1})
+c_{m-1}^2D(c_1,\ldots,c_{m-3})}{D(c_1,\ldots,c_{m-1})}\right)
-1 -\dt_{im}%
\\
&=&
 -2\frac{D(c_1,\ldots,c_{m-2})} {D(c_1,\ldots,c_{m-1})}
 -1 -\dt_{im}\\
 &=&
  -2J_{mm}-1 -\dt_{im},
  \eean
  confirming the last formula in the statement of Lemma 3.1.\qed

\begin{proposition}
For arbitrary $z \in \BC$, consider the map of
Proposition 2.2, namely

$$
((s_1,\ldots,s_m),(u^{(1)},\ldots,u^{(m)}))\longmapsto
((c_1,\ldots,c_{m-1},\al),(b^{(1)},\ldots,b^{(m)}))
$$
where

\bea  c_j&=& \left\{\begin{array}{ll}
\sqrt{\frac{(1+s_1z^2)(s_3-s_2)}{(1+s_2z^2)(s_3-s_1)}}
 &\mbox{for~}j=1\\
\\
\sqrt{\frac{(s_j-s_{j-1})(s_{j+2}-s_{j+1})}{(s_{j+1}-s_{j-1})
(s_{j+2}-s_j)}}&\mbox{for~}2\leq j\leq
m-2\\
\\
\sqrt{\frac{(s_{m-1}-s_{m-2})(1-s_mz^2)}{(s_m-s_{m-2})(1-s_{m-1}z^2)}}&\mbox{for~}j=m
-1
\end{array}\right.
 \no\\
 \no\\
\al
&=&\frac{2}{z}\sqrt{\frac{s_m-s_{m-1}}{(1-s_{m-1}z^2)(1-s_mz^2)}}
\no\\
\no\\
b^{(i)}&=&\frac{u^{(i)}}{U_i(s,z)}\qquad 1\leq i\leq m,
\label{F}\eea with

 $$
 U_i(s;z):=
 \left\{
  \begin{array}{ll }
 \displaystyle{\frac{1}{2}
\sqrt{\frac{(1+s_1z^2)(s_2-s_1)}{1+s_2z^2}}}
 ,&~~\mbox{~for~}
i=1,\\
\\
 \displaystyle{\frac{1}{2}
\sqrt{\frac{(s_i-s_{i-1})(s_{i+1}-s_i)}{s_{i+1}-s_{i-1}}}}
,&~~\mbox{~for~}
2\leq i\leq m-1,\\
\\
 \displaystyle{\frac{1}{2}\sqrt{\frac{(1-s_mz^2)(s_m-s_{m-1})}
 {1-s_{m-1}z^2}}}
,&~~\mbox{~for~}   i=m ,\end{array}
 \right.
  $$


The inverse map involves the tridiagonal matrix $J^{-1}$
and can be expressed as a fractional linear map in
$\frac{J_{m1}J_{mi}}{J_{1i}}$ for $1\leq i\leq m$
\bea
s_i&=&s_i(\al,c;z)=\frac{1}{z^2}\frac{\al^2z^4\frac{J_{m1}J_{mi}}{J_{1i}}
+2 }{\al^2z^4\frac{J_{m1}J_{mi}}{J_{1i}}-2 } \no\\\no \\
 u^{(i)}&=&b^{(i)}U_i(s(\al,c;z);z)
\label{G}
 \eea
 with
 $$
   U_i(s(\al,c;z);z)
 =  \frac{-\al z~J_{m1}}{\al^2z^4J_{m1}J_{mi}-2J_{1i}}
  =
  \frac{-1}{J_{mi}}\left(\frac{\al z~\frac{J_{m1}J_{mi}}{J_{1i}}}
 {\al^2z^4\frac{J_{m1}J_{mi}}{J_{1i}}-2}\right)
 $$
 Note that also the entries of $J$ can be expressed as
 functions of $s_i$:
  \bea
J_{1i}&=\displaystyle{\frac{-(1-z^2s_i)}{4z^2U_i(s,z)}\sqrt{\frac{(1+s_1z^2)(1+s_2z^2)}{s_2-s_1}}
}\no\\
J_{mi}&=\displaystyle{\frac{-(1+z^2s_i)}{4z^2U_i(s,z)}
\sqrt{\frac{(1-s_{m-1}z^2)(1-s_mz^2)}{s_m-s_{m-1}}}}
\label{H} \eea with $U_i(s,z)$ as in (\ref{F}).

\end{proposition}

 \proof  \underline {Step 1}:
 Upon expanding $D(c_1,\ldots,c_{m-1})$ along the
$j$-th row, one checks by (\ref{B}) and (\ref{C}), that
for $1\leq j\leq m-1$,
\bea
r_{j+1}-r_j&=&\frac{\displaystyle{\prod_1^{m-1}}(-c_i)}
{c_1^2\ldots c_j^2}\no\\
\no\\
&
&\frac{D(c_1,\ldots,c_{j-1})D(c_{j+1},\ldots,c_{m-1})-c_j^2D(c_{j+2},\ldots,c_{m-1})D(
c_1,\ldots,c_{j-2})}{D(c_{j+2},\ldots,c_{m-1})
D(c_{j+1},\ldots,c_{m-1})}  \no\\
\no\\
&=&
\frac{\displaystyle{\prod_1^{m-1}}(-c_i)}{\displaystyle{\prod_1^j}c_i^2}\cdot
\frac{D(c_1,\ldots,c_{m-1})}{D(c_{j+1},\ldots,c_{m-1})D(c_{j+2},\ldots,c_{m-1})},
 \label{r1}\eea  
and similarly, for $1\leq j\leq m-2$,
\bean
r_{j+2}-r_j&=&\frac{\displaystyle{\prod_1^{m-1}}(-c_j)}{\displaystyle{\prod_1^{j+1}}c_i^2}\\
\\
&
&\frac{D(c_1,\ldots,c_{j})D(c_{j+1},\ldots,c_{m-1})-c_j^2c^2_{j+1}D(c_{j+3},\ldots,c_{m-1})
D(
c_1,\ldots,c_{j-2})}{D(c_{j+3},\ldots,c_{m-1})D(c_{j+1},\ldots,c_{m-1})}\\
\\
&=&\frac{\displaystyle{\prod_1^{m-1}}(-c_i)}{\displaystyle{\prod_1^{j+1}}c_i^2}
\frac{D(c_1,\ldots,c_{m-1})}{D(c_{j+1},\ldots,c_{m-1})D(c_{j+3},\ldots,c_{m-1})}.
\eean
These identities then lead to:
$$
\frac{(r_i-r_{i-1})(r_{i+2}-r_{i+1})}{(r_{i+1}-r_{i-1})(r_{i+2}-r_i)}=c^2_i,\mbox{~for~}1\leq
i\leq m-2
$$
 $$
\frac{(r_i-r_{i-1})(r_{i+1}-r_{i})}{(r_{i+1}-r_{i-1})}=-\frac{J_{m1}}{J_{1i}^2}
 ,\mbox{~for~}1\leq
i\leq m-1
$$
and
\be \frac{r_{m-1}-r_{m-2}}{r_m-r_{m-2}}=c^2_{m-1}.
\label{E} \ee
$$
 r_m-r_{m-1}=-\frac{1}{J_{1m}}
  $$


\noindent \underline {Step 2}: It is easier to show that
the inverse map of (\ref{G}) is given by (\ref{F}). So,
from inverting the fractional linear map, appearing in
(\ref{G}), one computes
\be
J_{m1}r_i=\frac{J_{m1}J_{mi}}{J_{1i}}=-\frac{2}{\al^2z^4}
\left(\frac{1+z^2s_i}{1-z^2s_i}\right)\label{I} \ee and
so, inverting this map, one computes, in cascade,
\bean
1-z^2s_i&=&-\frac{4}{\al^2z^4r_iJ_{m1}-2}, \qquad 1+z^2s_i
 =\frac{2\al^2z^4r_iJ_{m1}}{\al^2z^4r_iJ_{m1}-2}\\
\\
s_i-s_{i-1}&=&\frac{4\al^2z^4(r_{i-1}-r_i)J_{m1}}
{(\al^2z^4r_iJ_{m1}-2)(\al^2z^4r_{i-1}J_{m1}-2)}\\
\\
s_{i+1}-s_{i-1}&=&\frac{4\al^2z^2(r_{i-1}-r_{i+1})
J_{m1}}{(\al^2z^4r_{i-1}J_{m1}-2)(\al^2z^4r_{i+1}J_{m1}+2)}.
\eean
 Therefore, using (\ref{E}), one checks
\bean
\frac{(s_i-s_{i-1})(s_{i+2}-s_{i+1})}{(s_{i+1}-s_{i-1})(s_{i+2}-s_i)}&=&
\frac{(r_i-r_{i-1})(r_{i+2}-r_{i+1})}{(r_{i+1}-r_{i-1})(r_{i+2}-r_i)}=c_i^2\mbox{~for~}2\leq
i\leq m-2\\
\\
\frac{(1+s_1z^2)(s_3-s_2)}{(1+s_2z^2)(s_3-s_1)}&=&\frac{r_1(r_3-r_2)}{r_2(r_3-r_1)}=c_1^2\\
\\
\frac{(1-s_mz^2)(s_{m-1}-s_{m-2})}{(1-s_{m-1}z^2)(1-s_mz^2)}&=&\frac{r_{m-1}-r_{m-2}}{r_m-r_{m-2}}=
c_{m-1}^2\\
\\
\frac{4}{z^2}\frac{s_m-s_{m-1}}{(1-s_{m-1}z^2)(1-s_mz^2)}
&=&-J_{m1}(r_m-r_{m-1})\al^2 =\al^2
\eean
and similarly for the expressions $U_i$ in (\ref {F}).
The signs are all specified from $\sqrt{s_{i+1}-s_i}$
and $\sqrt{ {s_{i+2}-s_i}}$.

Identity (\ref{H}) is obtained by solving 
$$
s_i=s_i(\al,c;z)=\frac{1}{z^2}
\frac{\al^2z^4\frac{J_{m1}J_{mi}}{J_{1i}}+2}
{\al^2z^4\frac{J_{m1}J_{mi}}{J_{1i}}-2 }
$$
for $\frac{J_{m1}J_{mi}}{J_{1i}}$, by substituting the
result in (\ref{G}), i.e.,
$$
 U_i(s(\al,c;z);z)
 = -\frac{ \al z~\frac{J_{m1}}{J_{1i}}  }
 {\al^2z^4\frac{J_{m1}J_{mi}}{J_{1i}}-2}
,$$
  and then solving for $\frac{J_{m1}}{J_{1i}}$.  Finally
  expressing
   $J_{mi}=\frac{\frac{J_{m1}J_{mi}}{J_{1i}}}
   {\frac{J_{m1}}{J_{1i}}}$
leads to the second relation (\ref{H}). The first
relation (\ref{H}) is obtained from the ratio
$J_{1i}=\frac{J_{m1}}{ {J_{m1}}/{J_{1i}}}$ and using the
result previously obtained for $J_{mi}$ at $i=1$. \qed


\begin{lemma}The following identities hold

\bean
 & &\al\frac{\pl}{\pl \al}s_i(\al,c;z)=
-8U_i^2(s,z)\frac{J_{mi}J_{1i}}{J_{m1}}
 = \frac{1}{z^2}-s_i^2z^2\\
\\
& &\al\frac{\pl}{\pl \al}\log U_i (s(\al,c;z),z)=-z^2s_i\\
\\
& &c_1\frac{\pl}{\pl c_1}s_i(\al,c;z)
 =8U_i^2(s,z)J^2_{1i}=\frac{(1-s_iz^2)^2(1+s_1z^2)
   (1+s_2z^2)}{2z^4(s_2-s_1)}\\
\\
& &c_1\frac{\pl}{\pl c_1}\!\log U_i(s(\al,c;z),z)
=\dt_{i1}-2\al z^3J_{m1}J_{1i}U_i(s,z)
\\
&&~~~~~~~~~~~~~~~~~~~~~~
~~~~~~~=\dt_{i1}\!-\!\frac{1}{2z^2(s_2\!-\!s_1)}(1+s_1z^2)(1+s_2z^2)(1-s_iz^2)
\\
& &c_{m-1}\frac{\pl}{\pl
c_{m-1}}s_i(\al,c;z)\\
&&~~~~=8U_i^2(s,z)\frac{J_{mi}J_{1i}}{J_{m1}}
\left(1+2J_{mm}-\frac{J_{mi}J_{m1}}
{J_{1i}}\right)\\
 &&
 ~~~~
  =
 -\frac{1+s_iz^2 }{2 z^4(s_m-s_{m-1})}\\
 &&~~~~~~
 \left(\left(1-{  s_m}z^{2}\right)\left(1-{ s_{m-1}}z^{2}\right)\,%
 \left({ 1+ s_i}z^{2} \right)-2
  \left({ 1-s_m s_{m-1}}z^{4} \right)
  \left({ 1- s_i}z^{2} \right)\right)
\\
& &c_{m-1}\frac{\pl}{\pl c_{m-1}}\log
U_i(s_i(\al,c;z),z)\\
&&~~~~=-1-2
\frac{D(c_{i+1},\ldots,c_{m-2})}{D(c_{i+1},\ldots,c_{m-1})}
 -2\al z^3U_i(s,z)
J_{mi}\left(1+2J_{mm}-\frac{J_{mi}J_{m1}}
{J_{1i}}\right)\\
 %
\eean
where
 \bean
 \lefteqn{ 1+2J_{mm}-\frac{J_{mi}J_{m1}}{J_{1i}}
 }\\&=&
 \frac{1}{2z^2(s_m-s_{m-1})(1-s_iz^2)}\\ \\
 &&
\times\left((1-s_mz^2)(1-s_{m-1}z^2)(1+s_iz^2)-
2(1-s_{m}s_{m-1}z^4)(1-s_iz^2)\right)
 \eean
 \bean
  \frac{D(c_{i+1},\ldots,c_{m-2})}{D(c_{i+1},\ldots,c_{m-1})}
  &=& \frac{J_{mm}J_{1i}-J_{m1}J_{mi}}{J_{1i}}
   \\&=&-\frac{(s_m-s_i)(1-z^2s_{m-1})}{(s_m-s_{m-1})(1-z^2s_{i})}
   \\
    J_{mm}&=&-\frac{(1-s_{m-1}z^2)(1+s_{m}z^2)}{2z^2(s_m-s_{m-1})}
\\
J_{m1}&=&-\frac{1}{2z^2}\sqrt{\frac{
{(1+s_1z^2)(1+s_{2}z^2)
  (1-s_{m-1}z^2)(1-s_{m}z^2)}}
   { {(s_2-s_1)(s_{m}-s_{m-1})}}}
   \\
 \frac{J_{mi}J_{m1}}{J_{1i}}&=& -{{\left(1-{ s_m} z^{2}\right) \left(1-{s_{m-1}}\,z^{2}\right)%
 \,\left(1+{  s_i} z^{2} \right)}\over{2z^{2}
  \left({  s_m}-{  s_{m-1}}%
 \right)\left(1-{  s_i} z^{2}\right)}}
  \\
J_{mi}U_i\al&=&\frac{-(1+s_iz^2)}{2z^3}
\\
J_{1i}U_i&=&-\frac{1-s_iz^2}{4z^2}
\sqrt{\frac{(1+s_1z^2)(1+s_2z^2)}{s_2-s_1}}
 .\eean

\end{lemma}


\proof Differentiating the first identity (\ref{G}) with
regard to $a$ and using the second expression (\ref{G})
yield
 $$
a\frac{\pl}{\pl a}s_i(a,c;z)= -{{8
a^{2}z^{2}\frac{J_{m1}J_{mi}}{J_{1i}} }\over{\left(a^{2}
 z^{4}\frac{J_{m1}J_{mi}}{J_{1i}} -2\right)^{2}}}
=-8U_i^2(s,z)\frac{J_{mi}J_{1i}}{J_{m1}}.
$$
Differentiating $U_i$ as in (\ref{G}), and using the
expression (\ref{G}) for $s_i$, yield:
 $$
 \al\frac{\pl}{\pl \al}\log U_i(\al,c;z)=
 -{\al^{2}{z^{4}\frac{J_{m1}J_{mi}}{J_{1i}}+2}
  \over
   {\al^{2}z^{4}\frac{J_{m1}J_{mi}}{J_{1i}}-2}}=-z^2 s_i
 ,$$
 while (\ref{G}) and Lemma 3.1 yield
 \bean
 c_1\frac{\pl}{\pl c_1}s_i(\al,c;z)
  &=&
   {{-4 \al^{2} z^{2%
 }}\over{\left(\al^{2} \frac{J_{m1}J_{mi}}{J_{1i}}\,z^{4}-2\right)^{2}}}
 c_1\frac{\pl }{\pl c_1}\frac{J_{m1}J_{mi}}{J_{1i}}
 \\
  &=&
   {{8 \al^{2} z^{2}J_{m1}^2}\over{\left(\al^{2}
\frac{J_{m1}J_{mi}}{J_{1i}}\,z^{4}-2\right)^{2}}}
\\
&=& 8U_i(s,z)^2J_{1i}^2\eean

 \bean
 c_{m-1}\frac{\pl}{\pl c_{m-1}}s_i(\al,c;z)
  &=&
   {{-4 \al^{2} z^{2%
 }}\over{\left(\al^{2} \frac{J_{m1}J_{mi}}{J_{1i}}\,z^{4}-2\right)^{2}}}
 c_{m-1}\frac{\pl}{\pl c_{m-1}}\frac{J_{m1}J_{mi}}{J_{1i}}
\\
&=&  8U_i(s,z)^2\frac{J_{mi}J_{1i}}{J_{m1}}
 \left(1+2J_{mm}-\frac{J_{m1}J_{mi}}{J_{1i}}\right) \eean

\bean
 c_1\frac{\pl}{\pl c_1}\log U_i(s;z)
  &=&
  c_1\frac{\pl}{\pl c_1}
   \log\frac{J_{m1} }{J_{1i}}-
   {{  \al^{2} z^{4} }\over{ \al^{2}
\frac{J_{m1}J_{mi}}{J_{1i}}\,z^{4}-2  }}
  c_1\frac{\pl}{\pl c_1}
   \frac{J_{m1}J_{mi}}{J_{1i}}
  \\
  &=&\dt_{1i}-2\al z^3J_{m1}J_{i1}U_i(s,z)
  \eean

  \bean
  \lefteqn{c_{m-1}\frac{\pl}{\pl c_{m-1}}\log
  U_i(s;z)}\\
  &=&
  c_{m-1}\frac{\pl}{\pl c_{m-1}}
   \log\frac{J_{m1} }{J_{1i}}-
   {{  \al^{2} z^{4} }\over{ \al^{2}
\frac{J_{m1}J_{mi}}{J_{1i}}\,z^{4}-2  }}
    c_{m-1}\frac{\pl}{\pl c_{m-1}}
   \frac{J_{m1}J_{mi}}{J_{1i}}
  \\
  &=&-1-2
\frac{D(c_{i+1},\ldots,c_{m-2})}{D(c_{i+1},\ldots,c_{m-1})}\\
\\& &\hspace*{2cm}-2\al z^3U_i(s,z)
    J_{mi}  \left(1+2J_{mm}-\frac{J_{mi}J_{m1}}
{J_{1i}}\right)\eean yielding most of the differential
identities of Lemma 3.3. The remaining relations are a
consequence of (\ref{F}), (\ref{G}), (\ref{H}) and
(\ref{I}). \qed

\newpage

\section{Integrable deformations and the Virasoro constraints}

In order to compute the differential equation for
\bea
\lefteqn{\BP_n\left(\alpha;c_1,\ldots,c_{m-1};E_1\times
\ldots\times E_m\right)
}\no\\
& =& \frac{1}{Z_n}\int_{\prod_1^m E_i^{k_1+k_2}}
\Dt_{k_1+k_2}(x^{(1)})
   \prod_{i=1}^{k_1+k_2} \prod_{\ell=1}^{m}
   dx_i^{(\ell)}
 \no\\
   && \Dt_{k_1}(x^{(m)'})
\prod_{i=1}^{k_1}e^{-\frac{1}{2}\sum_{\ell=1}^{m}
x_i^{(\ell)2}
  + \sum_{\ell=1}^{m-1}c_{\ell}x_i^{(\ell)}x_i^{(\ell+1)}
   +\alpha x_i^{(m)}  }\no\\
 &&   \Dt_{k_2}(x^{(m)''})
  \!\!\prod_{i=k_1+1}^{k_1+k_2}\!\!
   e^{-\frac{1}{2}\sum_{\ell=1}^{m}
x_i^{(\ell)2}
  + \sum_{\ell=1}^{m-1}c_{\ell}x_i^{(\ell)}x_i^{(\ell+1)}
   -\alpha x_i^{(m)}},
  \label{4.A} \eea
   we need to add to the numerator of $\BP_n$ many auxiliary variables
   $$
  \left\{ \begin{array}{l}
    \bar t:=(\bar t_1,\bar t_2,\ldots),\bar s:=(\bar s_1,\bar s_2,\ldots),
    \bar u:=(\bar u_1,\bar u_2,\ldots)\mbox{  and  }
    \beta
 ,\\ \\
 \gamma^{(\ell)}:=(\gamma_1^{(\ell)},\gamma_2^{(\ell)},\ldots)
 \mbox{  for } 2\leq \ell\leq m-1
 \\ \\
  c^{(\ell)}:=(c_{p,q}^{(\ell)})_{p,q\geq 1}  \mbox{  for } 1\leq \ell\leq m-1$$
\end{array}\right.
.$$ Note that the time variables $\bar t,~\bar s,~\bar
u$ are totally different from the $t$-variables
appearing in the Brownian motion.
 yielding the following integral\footnote{If $m=1$ or
 $2$, the formulae below must be reinterpreted; e.g.,
 for $m=1$, the $c_i$'s are absent and for $m=2$, $\gamma$
 is not present.}
\bea \lefteqn{\tau_{k_1k_2}(\bar t,\bar s,\bar
u;\beta,\gamma^{(2)},\ldots, \gamma^{(m-1)},c
^{(1)},\dots,c ^{(m-1)}, \alpha, E_1\times
\ldots\times E_m)}\no\\
 \no\\
 &=&
  \frac{1}{k_1!k_2!}\int_{\prod_1^m
E_i^{k_1+k_2}} \Dt_{k_1+k_2}(x^{(1)})
 \prod_{i=1}^{k_1+k_2}\left(e^{\sum_{j=1}^{\iy}
\bar
t_{j}x_i^{(1)j}}\prod_{\ell=1}^{m}dx_i^{(\ell)}\right)
\no\\
\no\\
&&\Dt_{k_1}(x_1^{(m)},\ldots,x_{k_1}^{(m)})
 \no\\\no\\
&&\prod_{i=1}^{k_1}e^{-\frac{1}{2}\sum_{\ell=1}^{m}
x_i^{(\ell)2}+\alpha x_i^{(m)}+\beta
x_i^{(m)2}-\sum_{j=1}^{\iy}\bar s_{j}(x_i^{(m)})^{j}
 +\sum_{ p,q\geq 1} \sum_{\ell=1}^{m-1}c^{(\ell)}_{ p,q}
 (x_i^{({\ell})})^{p}
 (x_i^{({\ell}+1)})^{q}
   +\sum_{{\ell}=2}^{m-1}\sum_{r=1}^{\iy}\gamma^{(\ell)}_r
     (x_i^{({\ell})})^r
     }
   \no\\\no\\
  &&
  \Dt_{k_2}(x^{(m)}_{k_1+1},\ldots,x^{(m)}_{k_1+k_2})
  \no\\\no\\
  &&
 \!\!\!\!\prod_{i=k_1+1}^{k_1+k_2}\!\!\!
  e^{-\frac{1}{2}\sum_{\ell=1}^{m}
x_i^{(\ell)2}-\alpha x_i^{(m)}-\beta
x_i^{(m)2}-\sum_{j=1}^{\iy}\bar u_{j}(x_i^{(m)})^{j}
 +\sum_{ p,q\geq 1} \sum_{\ell=1}^{m-1}c^{(\ell)}_{ p,q}
 (x_i^{({\ell})})^{p}
 (x_i^{({\ell}+1)})^{q}
   +\sum_{{\ell}=2}^{m-1}\sum_{r=1}^{\iy}\gamma^{(\ell)}_r
     (x_i^{({\ell})})^r}
   \no\\
    &=& \det\left(  \begin{array}{l}
        (\mu_{ij}^+)_{1\leq i \leq k_1,~1\leq j\leq
        k_1+k_2 }
        \\  \\
        (\mu_{ij}^-)_{1\leq i \leq k_2,~1\leq j\leq
        k_1+k_2 }
        \\
        \end{array}
        \right)
\label{4.B}\eea
where
  \bea
   \mu_{ij}^{\pm}
   &=&
   \int_{\prod_1^m E_i}
  \left( \prod_{\ell=1}^m dx^{(\ell)}\right)
   e^{\sum_1^{\iy}\left(\bar t_kx^{(1)k}
    -\left({{\bar s_k}\atop{ \bar u_k}}\right)
    x^{(m)k}\right)}
 (x^{(1)})^{j-1} (x^{(m)})^{i-1}
  \no\\\no\\
 && ~~ e^{-\frac{1}{2}\sum_{\ell=1}^{m} x^{(\ell)2}
 \pm\alpha x^{(m)}\pm\beta x^{(m)2}+ \sum_{ p,q\geq 1}\sum_{\ell=1}^{m-1}c^{(\ell)}_{pq}(x^{(\ell)})^p
 (x^{(\ell+1)})^q
  +
  \sum_{{\ell}=2}^{m-1}\sum_{r=1}^{\iy}\gamma^{(\ell)}_r
     (x^{({\ell})})^r
   }
  \no\\
  \label{4.C}\eea
This is to say, the integral (\ref{4.B}), along the
locus
\be
\LR:=\left\{\begin{array}{l}\bar t_i=0,~\bar s_i=0,\bar u_i=0,
\beta=0,~\gamma_r^{(\ell)}=0,\\
 c_{11}^{(\ell)}=c_{\ell}~~
\mbox{and}~~c_{ij}^{(\ell)}=0~~\mbox{for}~~i,j\geq  1
~\mbox{with}~ (i,j)\neq (1,1)
\end{array}\right\}
\label{locus}\ee
  yields the integral (\ref{4.A}); also for the sake
of brevity, set $\gamma_{\ell}:=\ga_1^{(\ell)}$.
The following locus \be
\LR_{\beta}:=\left\{\begin{array}{l}\bar t_i=0,~\bar s_i=0,
\bar u_i=0,~\gamma_r^{(\ell)}=0,\\
 c_{11}^{(\ell)}=c_{\ell}~~
\mbox{and}~~c_{ij}^{(\ell)}=0~~\mbox{for}~~i,j\geq  1
~\mbox{with}~ (i,j)\neq (1,1)
\end{array}\right\}
\label{locusbeta}\ee will also be used. Then the
following statement holds:

\begin{proposition}
  Given a disjoint union of intervals and the associated
  algebra of differential operators
  $$
  E_{\ell}:=
\bigcup^r_{i=1}[b^{(\ell)}_{2i-1},b^{(\ell)}_{2i}]\subset
\BR ~~\mbox{and}~~
 {\cal
   D}_k(E_{\ell})=\sum_{i=1}^{2r}(b_i^{(\ell)})^{k+1}\frac{\pl}{\pl
b_i^{(\ell)}},\qquad 1\leq\ell\leq m$$
 the integrals (\ref{4.B})
satisfy, besides
 the trivial
relations,
\be
  -\frac{\pl}{\pl \bar s_1}+\frac{\pl}{\pl
\bar u_1}=\frac{\pl}{\pl \alpha}
 ,\qquad -\frac{\pl}{\pl \bar s_2}+\frac{\pl}{\pl
\bar u_2}=\frac{\pl}{\pl \beta} ,
 \label{trivial}\ee
 and, upon setting
 $$\gamma_1:=\bar t_1,\mbox{   and} ~~~\frac{\pl}{\pl
 \gamma_m}:=
  -\frac{\pl}{\pl
\bar s_1}-\frac{\pl}{\pl \bar u_1} ,$$
 the following Virasoro relations
 ($2\leq \ell\leq m-1$):
\bean
\DR_{-1}(E_{1})\tau_{k_1k_2}&=&\left\{\begin{array}{l}
 \displaystyle{\stackrel{ -\frac{\pl}{\pl
\gamma_1}}{\overbrace{-\frac{\pl}{\pl \bar t_1}}}+c_1
\frac{\pl}{\pl
\gamma_2} +(k_1+k_2)\bar t_1 }  \\ \\
 \displaystyle{+\sum_{i\geq 2} i\bar t_i\frac{\pl}{\pl
\bar t_{i-1}}+
  \sum_{{i\geq 2}\atop{j\geq 1}}ic_{ij}^{(1)}
 \frac{\pl}{\pl c^{(1)}_{i-1,j}}
  + \sum_{{j\geq 2}} c_{1j}^{(1)}
 \frac{\pl}{\pl \gamma^{(2)}_{j}}
 }
\end{array}\right\}\tau_{k_1k_2}
\\
\vdots
 \eean
 \bea
\DR_{-1}(E_{\ell})\tau_{k_1k_2}&=&\left\{\begin{array}{l}
\displaystyle{c_{\ell-1}\frac{\pl}{\pl
\gamma_{\ell-1}}-\frac{\pl}{\pl
\gamma_{\ell}}+c_{\ell}\frac{\pl}{\pl
\gamma_{\ell+1}}+(k_1+k_2)\gamma_{\ell}
}\no\\\no\\
\displaystyle{
+ \sum_{{i\geq 1}\atop{j\geq
2}}jc^{(\ell-1)}_{ij}\frac{\pl}{\pl
 c^{(\ell-1)}_{i,j-1}} +
  \sum_{{i\geq 2}\atop{j\geq 1}}ic^{(\ell )}_{ij}
  \frac{\pl}{\pl c^{(\ell )}_{i-1,j}}
 }\\ \\
 +\sum_{{i\geq 2}} c^{(\ell-1)}_{i1}\frac{\pl}{\pl
 \gamma^{(\ell-1)}_{i}}
  +\sum_{{j\geq 2}}  c^{(\ell)}_{1j}\frac{\pl}{\pl
 \gamma^{(\ell+1)}_{j}}
  \\ \\
  +\sum_{{r\geq 2}}r \gamma^{(\ell)}_{r}\frac{\pl}{\pl
 \gamma^{(\ell)}_{r-1}}
\end{array}\right\}\tau_{k_1k_2}
\no\\\no\\
\vdots&&\no\\
\no\\
\DR_{-1}(E_{m})\tau_{k_1k_2}&=&\left\{\begin{array}{l}
 \displaystyle{c_{m-1}\frac{\pl}{\pl \gamma_{m-1}}
+ \stackrel{ -\frac{\pl}{\pl
\gamma_m}}{\bigl(\overbrace{\frac{\pl}{\pl \bar
s_1}\!+\!\frac{\pl}{\pl \bar u_1}}\bigr)}+2\beta
\stackrel{\frac{\pl}{\pl \al}}{
(\overbrace{-\frac{\pl}{\pl \bar s_1}\!+\!\frac{\pl}{\pl
\bar u_1}})}}
\\ \\
 \displaystyle{-k_1(\bar s_1-\alpha)-k_2(\bar u_1+\alpha)}
  \\ \\
\displaystyle{+\sum_{i\geq 2}i\left( \bar
s_i\frac{\pl}{\pl \bar s_{i-1}}+ \bar u_i\frac{\pl}{\pl
\bar u_{i-1}}\right)
 }\\+\displaystyle{\sum_{{i\geq 2}} c^{(m-1)}_{i1}\frac{\pl}{\pl
 \gamma^{(m-1)}_{i}}
  +
\sum_{{i\geq 1}\atop{j\geq
2}}jc^{(m-1)}_{ij}\frac{\pl}{\pl c^{(m-1)}_{i,j-1}}}
\end{array}\right\}\tau_{k_1k_2}\no\\
\label{4.Vir1}\eea
and (only needed for $\ell=1,m$)
 \bean
\DR_0(E_{1})\tau_{k_1k_2}&=&\left(
 \begin{array}{l}
 \displaystyle{-\frac{\pl}{\pl \bar t_2} +c_1\frac{\pl}{\pl
c_1} + \frac{(k_1+k_2)(k_1+k_2+1)}{2}}\\  \\
\displaystyle{\sum_{i\geq 1} i\bar t_i\frac{\pl}{\pl
\bar t_{i}}+ \sum_{{i,j\geq 1}\atop{ (i,j)\neq (1,1) }
}ic^{(1)}_{ij}\frac{\pl}{\pl c^{(1)}_{i,j}}}
\end{array}
\right) \tau_{k_1k_2}
\eean
 \bea
  \DR_0(E_{m})\tau_{k_1k_2}&=&\left(\begin{array}{l}
\displaystyle{ \alpha\bigl(\stackrel{\frac{\pl}{\pl
\al}}{\overbrace{-\frac{\pl}{\pl \bar
s_1}+\frac{\pl}{\pl \bar u_1}}} \bigr)
 +c_{m-1}\frac{\pl}{\pl c_{m-1}}
  +\bigl(\frac{\pl}{\pl \bar s_2}+ \frac{\pl}{\pl \bar u_2}\bigr)}
 \no\\ \no\\
 \displaystyle{
+2\beta \bigl(\stackrel{\frac{\pl}{\pl
\beta}}{\overbrace{-\frac{\pl}{\pl \bar s_2}+
\frac{\pl}{\pl
\bar u_2}}}\bigr)+\frac{k_1(k_1+1)}{2}+\frac{k_2(k_2+1)}{2}}\\ \\
\displaystyle{\sum_{i\geq 1}\left(i\bar
s_i\frac{\pl}{\pl \bar s_{i}}+i\bar u_i\frac{\pl}{\pl
\bar u_{i}}\right)+ \sum_{{i,j\geq 1}\atop{(i,j)\neq
(1,1)}}jc^{(m-1)}_{ij}\frac{\pl}{\pl c^{(m-1)}_{i,j}}}
\end{array}\right)\tau_{k_1k_2}\no\\
\no\\\label{4.Vir2}
 \eea

 \end{proposition}


Before giving the proof of Proposition 4.1, we need the
following lemma, concerning the expressions:

$$
dI_n:=\Dt_n(x)\prod_{k=1}^n dx_k e^{\sum_{i=1}^{\iy}\bar
t_ix_k^i}
$$
$$
 {\cal E}(x,y):= \prod_{k=1}^n
e^{\sum_1^{\iy}\bar t_ix_k^i+\sum_{i,j \geq 1}c_{ij}
x_k^i y_k^j+\sum_1^{\iy}\bar s_iy_k^i}
$$

\begin{lemma}

$$
 \frac{\pl}{\pl \vr}
  \prod_{i=1}^n d(
x_i +\vr x_i^{k+1})
   \Bigr|_{\vr =0} =\left\{
 \begin{array}{ll}
 0,&k=-1\\
 n\prod_{i=1}^n dx_i,&k=0
 \end{array}\right.
 $$

$$
 \frac{\pl}{\pl \vr}
 dI_n(
x +\vr x^{k+1})\Bigr|_{\vr =0} =\left\{
 \begin{array}{l}\left(\sum_{i\geq 2} i\bar t_i\frac{\pl}{\pl
\bar t_{i-1}} +n\bar t_1\right)dI_n(x),~~~k=-1\\  \\
 \left(\sum_{i\geq 1} i\bar t_i\frac{\pl}{\pl
\bar t_{i }}+\frac{n(n+1)}{2}\right)dI_n(x),~~~k=0
\end{array}\right.
 $$

\bean \lefteqn{ \frac{\pl}{\pl \vr}
 {\cal E}(
x +\vr x^{k+1},y)\Bigr|_{\vr =0} }\\
&=& \left\{\begin{array}{ll}
\left(\displaystyle{\sum_{i\geq 2} i\bar
t_i\frac{\pl}{\pl \bar t_{i-1}}}+n\bar t_1+
  \sum_{i\geq 2,j\geq 1} ic_{ij}
\frac{\pl}{\pl c_{i-1,j}} +\sum_{j\geq 1}  c_{1j}
 \frac{\pl}{\pl \bar s_{j}}\right)
  {\cal E}(x,y) ,&k=-1\\
 \\
\left(\displaystyle{\sum_{i\geq 1} it_i\frac{\pl}{\pl
\bar t_{i }}}+
  \sum_{i,j\geq 1} ic_{ij}
\frac{\pl}{\pl c_{ij}} \right)
  {\cal E}(x,y),&k=0
  \end{array}
  \right.
 \eean

  \bean \lefteqn{
\frac{\pl}{\pl \vr}
 {\cal E}(
x, y +\vr y^{k+1})\Bigr|_{\vr =0} }\\
 &=&
\left\{\begin{array}{ll}
\left(\displaystyle{\sum_{i\geq 2} i\bar
s_i\frac{\pl}{\pl \bar s_{i-1}}}+n\bar s_1+
  \sum_{i\geq 1,j\geq 2} jc_{ij}
\frac{\pl}{\pl c_{i,j-1}} +\sum_{i\geq 1}  c_{i1}
 \frac{\pl}{\pl \bar t_{i}}\right)
  {\cal E}(x,y),&k=-1\\
 \\
\left(\displaystyle{\sum_{i\geq 1} i\bar
s_i\frac{\pl}{\pl \bar s_{i }}}+
  \sum_{i,j\geq 1} jc_{ij}
\frac{\pl}{\pl c_{ij}} \right)
  {\cal E}(x,y),&k=0
  \end{array}
  \right.
  \eean

  \end{lemma}

%


\proof This is obtained by setting
 $x_i\mapsto
x_i+\vr x_i^{k+1}$ and then $y_i\mapsto y_i+\vr
y_i^{k+1}$; then take the derivative with regard to
$\vr$ and set $\vr=0$.



\medskip\noindent
{\it Proof of Proposition 4.1:\/}

\noindent\underline{Case 1}:~Performing the
infinitesimal change of variables, for all $1\leq i\leq
n$,
 $$
 x^{(1)}_i\longmapsto x^{(1)}_i+\vr {x^{(1)}_i}^{k+1}
 $$
 in the integral (\ref{4.B}) involves the following integral only
$$\prod_{i=1}^{k_1+k_2}dx^{(1)}_i\Dt_{k_1+k_2}(x^{(1)})
  \prod_{i=1}^{k_1+k_2}e^{\sum_{j=1}^{\iy} \bar t_{j}(x_i^{(1)})^{j}}
\prod^{k_1}_{i=1}e^{-\frac{1}{2}x_i^{(1)2}+\sum_{p,q\geq
1}c_{p,q}^{(1)}(x_i^{(1)})^p(x_i^{(2)})^q +
\sum_{r=1}^{\iy}\gamma^{(2)}_r
     (x_i^{({2})})^r}
$$
and leads to the first Virasoro constraints in
(\ref{4.Vir1}), upon differentiation with regard to
$\vr$, setting $\vr =0$, applying Lemma 4.2 and taking
into account the variation of the boundary term in the
integral, $b_j^{(1)}\mapsto b_j^{(1)}-\vr
(b_j^{(1)})^{k+1}+O(\vr^2)$. The shifts in the time
parameter produces the terms in the first line of the
Virasoro constraint.

\noindent\underline{Case 2}:~Performing the
infinitesimal change of variables, for all $1\leq i\leq
n$,
 $$
 x_i^{(\ell)}\longmapsto x_i^{(\ell)}+\vr {x_i^{(\ell)}}^{k+1}
 $$
 in the integral (\ref{4.B}) involves the following integral
 only
\bean
\lefteqn{\prod_{i=1}^{k_1+k_2}dx^{(\ell)}_i
 \prod_{i=1}^{k_1}e^{ \sum_{r=1}^{\iy}\gamma^{(\ell)}_r
     (x_i^{({\ell})})^r-\frac{1}{2}x_i^{(\ell
)2}
+ \sum_{r=1}^{\iy}\gamma^{(\ell-1)}_r
     (x_i^{({\ell-1})})^r+ \sum_{p,q\geq 1}
 c_{p,q}^{(\ell-1)}(x_i^{(\ell-1)})^p(x_i^{(\ell)})^q
 }}\hspace{5cm}\\
\\
& &e^{  \sum_{r=1}^{\iy}\gamma^{(\ell+1)}_r
     (x_i^{({\ell+1})})^r+
 \sum_{p,q\geq 1}c^{(\ell)}_{p,q}
 (x_i^{(\ell)})^p (x_i^{(\ell+1)})^q}
\eean and leads to the Virasoro constraints for $2\leq
\ell\leq m-1$ in (\ref{4.Vir1}), upon differentiation
with regard to $\vr$, setting $\vr =0$, applying Lemma
4.2 and taking into account the variation of the
boundary term in the integral, $b_j^{(\ell)}\mapsto
b_j^{(\ell)}-\vr (b_j^{(\ell)})^{k+1}+O(\vr^2)$. Also
here one must take into account the shifts.

\noindent\underline{Case 3}:~Performing the
infinitesimal change of variables, for all $1\leq i\leq
n$,
 $$
 x_i^{(m)}\longmapsto x_i^{(m)}+\vr {x_i^{(m)}}^{k+1}
 $$
 in the integral (\ref{4.B}) involves the following integral
 only
%
%
%
 %
\bean \lefteqn{\prod_{i=1}^{k_1+k_2}dx^{(m)}_i
 \Dt_{k_1}(x_1^{(m)},\ldots,x_{k_1}^{(m)})
  \Dt_{k_2}(x^{(m)}_{k_1+1},\ldots,x^{(m)}_{k_1+k_2})}
 \\
\\
& &\hspace{-.7cm}
\prod_{i=1}^{k_1}e^{
- \sum^{\iy}_{j=1}
[\bar s_j-\dt_{j1}\al+\dt_{j2}(\frac{1}{2}-\beta)](x_i^{(m)})^j
+  \sum_{r=1}^{\iy} \gamma^{(m-1)}_r
     (x_i^{({m-1})})^r+\sum_{p,q\geq
1}c_{p,q}^{(m-1)}(x_i^{(m-1)})^p(x_i^{(m)})^q}\\
\\
& &\!\!\!
 \hspace{-.7cm}\prod_{i=k_1}^{k_1+k_2}
e^{
-\sum^{\iy}_{j=1}
 [\bar u_j+ \dt_{j1}\al+\dt_{j2}(\frac{1}{2}+\beta)] (x_i^{(m)})^j
+ \sum_{r=1}^{\iy}\gamma^{(m-1)}_r
     (x_i^{({m-1})})^r+
  \sum_{p,q\geq
1}c_{p,q}^{(m-1)}(x_i^{(m-1)})^p(x_i^{(m)})^q},
 \eean
%
%
and leads to the last Virasoro constraint for $\ell=m$
with a similar argument; this ends the proof of
Proposition 4.1. The relations (\ref{trivial}) follow at
once by inspection of the integral (\ref{4.B}).\qed

\newpage

Setting
\be
 J^{-1}:=
J^{-1}(c_1,\ldots,c_{m-1}):=\left(\begin{array}{ccccccccccc}
-1& &c_1\\
& & & \ddots& & {\bf O}\\
c_1& &-1& & \ddots& \\
 &\ddots& &\ddots & &\ddots\\
& & \ddots&  & -1& &c_{m-1}\\
 &{\bf O}& &\ddots\\
& & &  &  c_{m-1}& &-1
\end{array}
\right), \label{tridiag}
 \ee
 define the differential operators $
{\cal A}_i^{\pm},   {\cal C}_i$ and $\bar {\cal A}_i,
\bar {\cal C}_i$:
  \bean
  \!\!\AR_1^{\pm} &:=&-\frac{1}{2}
  \sum_{j=1}^m J_{mj} \left( \DR_{-1}(b^{(j)})
  -2\dt_{jm}\beta\frac{\pl}{\pl \alpha} \right)\mp
   \frac{1}{2}\frac{\pl}{\pl \alpha}
    =: \bar\AR_1^{\pm}
 +\beta J_{mm}\frac{\pl}{\pl \al}
  \\
  \CR_1&:=&\sum_{j=1}^m J_{1j} \left( \DR_{-1}(b^{(j)})
  -2\dt_{jm}\beta\frac{\pl}{\pl \alpha} \right)
  =:\bar\CR_1-2\beta J_{1m}\frac{\pl}{\pl \alpha}
  \\
  \\
  \\
  \!\!\AR_2^{\pm}
  \!\! &:=&\frac{1}{2}
 \left(
  \DR_0(b^{(m)})- \alpha\frac{\pl}{\pl \alpha} - c_{m-1}\frac{\pl}{\pl c_{m-1}}
  \right)- ( \beta \pm \frac{1}{2})\frac{\pl}{\pl \beta }
  =:  \bar\AR_2
 \mp\frac{1}{2}\frac{\pl}{\pl \beta
 }-\beta \frac{\pl}{\pl \beta
 }
 \\
{\cal C}_2&:=&-
  \DR_0(b^{(1)}) + c_1\frac{\pl}{\pl c_1}=:\bar {\cal C}_2
 \eean
We now state:

\begin{proposition} Along the locus $\LR$, the partials
and second partials of $f:=\log
 \tau_{k_1k_2}$ with regard to $\bar t_i,\bar s_i,\bar u_i$ can be
 expressed in terms of the operators
  $
\bar{\cal A}_i,   \bar{\cal C}_i$ and $\pl/\pl \beta$:

  \bea \bar{\cal A}_1^+f&=&\frac{\pl
f}{\pl \bar s_1} + \frac{\alpha \left({k_2}-{
k_1}\right)}{2}{{ } }J_{mm}
\no\\
\no\\
\bar{\cal A}_1^-f&=&\frac{\pl f}{\pl \bar u_1} +
{{\alpha\,\left({k_2}-{ k_1}\right)}\over{2}}
J_{mm}
\no\\
\no\\
\bar{\cal C}_1f&=&\frac{\pl f}{\pl \bar t_1}
 -\alpha (k_2-k_1)J_{m1}
 \no\\
\no\\
\left(\bar{\cal A}_2-\frac{1}{2}\frac{\pl}{\pl
\beta}\right)f&=&\frac{\pl f}{\pl \bar s_2} +
\frac{1}{4}({{ k_2}^{2}+{ k_2}+{ k_1}^{2}+{ k_1}})
\no\\
\no\\
\left(\bar{\cal A}_2+\frac{1}{2}\frac{\pl}{\pl
\beta}\right)f&=&\frac{\pl f}{\pl \bar u_2} +
\frac{1}{4}({{ k_2}^{2}+{ k_2}+{  k_1}^{2}+{  k_1}})
\no\\
\no\\
\bar{\cal C}_2f&=&\frac{\pl f}{\pl \bar t_2}
-\frac{1}{2}{ \left({ k_2}+{  k_1}\right)
\left({  k_2}+{  k_1}+1%
 \right)}
 \no\\
\no\\
  \bar{\cal A}_1^+\bar{\cal C}_1 f&=&\frac{\pl^2 f}
  {\pl \bar t_1\pl \bar s_1}
 -k_1 J_{1m} \no\\
\no\\
\left(\bar{\AR}_2\bar{\CR}_1+J_{1m} \frac{\pl}{\pl\al}
-\frac{1}{2}\bar{\CR}_1\frac{\pl  }{\pl\beta}\right)
f&=&\frac{\pl^2 f}{\pl \bar t_1\pl \bar s_2} -\alpha
(k_2-k_1) J_{1m}J_{mm}
\no\\
\no\\
\bar{\cal C}_2\bar{\cal A}_1^+ f&=&\frac{\pl^2 f}{\pl
\bar t_2\pl \bar s_1}
 -\alpha (k_2-k_1) J_{1m}^2
\label{Vir-identities}\eea

\end{proposition}

\proof The Virasoro relations (\ref{4.Vir1}) can be
written as (remembering $t_1=\gamma_1$)

$$\left(\begin{array}{c}
 {\cal D}_{-1}(b^{(1)})\\ \\
 \vdots
 \\ \\
 {\cal D}_{-1}(b^{(m-1)})\\  \\
{\cal D}_{-1}(b^{(m)})-2\beta\frac{\pl}{\pl \alpha}
 -(k_1-k_2)\al\end{array}\right)f
 =J^{-1}
  \left(\begin{array}{c}
 \frac{\pl}{\pl \gamma_{1}}\\ \\
 \vdots
 \\ \\
\frac{\pl}{\pl \gamma_{m}}
 \end{array}\right)f+O(\LR_\beta
 ),$$
  remembering the locus $\LR_\beta$ as in
  (\ref{locusbeta}),
and the tridiagonal matrix $J^{-1}$ defined in
(\ref{tridiag}). The symbol $O(\LR_\beta)$ means a term
which vanishes along $\LR_\beta$. Therefore, using the
trivial relations (\ref{trivial}), one finds
$$
J
 \left(\begin{array}{c}
 {\cal D}_{-1}(b^{(1)})\\ \\
 {\cal D}_{-1}(b^{(2)})\\ \\
 \vdots
 \\ \\
 {\cal D}_{-1}(b^{(m-1)})\\  \\
{\cal D}_{-1}(b^{(m)})-2\beta\frac{\pl}{\pl \alpha}
 -(k_1-k_2)\al\end{array}\right)
 f=
  \left(\begin{array}{c}
 \frac{\pl}{\pl t_{1}}\\ \\
 \frac{\pl}{\pl \gamma_{2}} \\ \\
 \vdots
 \\ \\
 \frac{\pl}{\pl \gamma_{m-1}}
 \\ \\
-\frac{\pl}{\pl \bar s_1}-\frac{\pl}{\pl \bar u_1}
 \end{array}\right)f+O(\LR_\beta)$$
 which combined with
  $$-\frac{\pl}{\pl \bar s_1}+\frac{\pl}{\pl
\bar u_1}= \frac{\pl}{\pl \alpha},$$
 leads to the following equations for
 $f=\log \tau_{k_1k_2}$ along the locus $\LR_\beta$:
%
$$
{\cal A}_1^+f=\frac{\pl f}{\pl \bar s_1} +
\frac{\alpha}{2}{{
\left({k_2}-{ k_1}\right)} }J_{mm}
$$

$$
{\cal A}_1^-f=\frac{\pl f}{\pl \bar u_1} + {{\alpha
}\over{2}}\left({k_2}-{ k_1}\right) J_{mm}
$$

$$
{\cal C}_1f=\frac{\pl f}{\pl \bar t_1}
 -\alpha (k_2-k_1)J_{m1}
 $$
The second set of Virasoro relations, combined with
$$-\frac{\pl}{\pl \bar s_2}+\frac{\pl}{\pl
\bar u_2}= \frac{\pl}{\pl \beta},$$
 leads to the following equations, also along the locus $\LR_\beta$,
\bean {\cal A}_2^+f&=&\frac{\pl f}{\pl \bar s_2} +
\frac{1}{4}({{ k_2}^{2}+{ k_2}+{  k_1}^{2}+{  k_1}})
\\
{\cal A}_2^-f&=&\frac{\pl f}{\pl \bar u_2} +
\frac{1}{4}({{ k_2}^{2}+{ k_2}+{  k_1}^{2}+{  k_1}})
\\
{\cal C}_2f&=&\frac{\pl f}{\pl \bar t_2}  -\frac{1}{2}{
\left({ k_2}+{  k_1}\right)
\left({  k_2}+{  k_1}+1%
 \right)}  \eean
Moreover
$$
{\cal A}_1^+{\cal C}_1 f=\frac{\pl^2 f}{\pl \bar t_1\pl
\bar s_1}
 -k_1 J_{1m} -\beta J_{mm}J_{m1}(k_2-k_1)$$


$$
{\cal A}_2^+{\cal C}_1 f=\frac{\pl^2 f}{\pl \bar t_1\pl
\bar s_2} -\alpha (k_2-k_1) J_{1m}J_{mm}
$$

$$
{\cal C}_2{\cal A}_1^+ f=\frac{\pl^2 f}{\pl \bar t_2\pl
\bar s_1}
 -\alpha (k_2-k_1) J_{1m}^2 . $$
Indeed, one first needs to check
$$
 {\cal A}_1 (\al)=J_{mm}\beta-\frac{1}{2}
 ,$$
 \bean
 {\cal A}_2^+(\al J_{m1})&=&
  -\frac{1}{2} \left( \al \frac{\pl}{\pl \al}+c_{m-1}\frac{\pl}{\pl c_{m-1}}
  \right)\al J_{m1}
  \\
  &=&
   -\frac{1}{2} \left(\al J_{m1} +\al J_{m1}(-2J_{mm}-1)\right)
\\
&=& \al J_{m1}J_{mm}
 \eean
and
 $$
  {\cal C}_2( J_{mm})=c_1\frac{\pl }{\pl c_1} J_{mm}=
  -2J_{m1}^2.$$
Then, using the expressions for the ${\cal A}_i^\pm f$
and ${\cal C}_i f$, one computes along the
locus\footnote{Since one takes the derivative with
regard to $\bar t_1$ in the expression below, one must
keep track of the term containing $\bar t_1$ in the
identity for ${\cal A}_1^+f$ before setting $\bar
t_1=0$.} $\LR_{\beta}$,
 \bean {\cal A}_1^+{\cal C}_1f&=&
 {\cal A}_1^+\left( \frac{\pl f}{\pl \bar t_1}
 -\alpha (k_2-k_1)J_{m1}\right)
 \\&=&\frac{\pl  }{\pl
\bar t_1}({\cal A}_1^+f)
 -  (k_2-k_1)J_{m1}({\cal A}_1^+\alpha)
\\
 &=&
  \frac{\pl  }{\pl \bar t_1} \left(\frac{\pl f}{\pl \bar s_1}
  + \frac{\alpha}{2}{{
\left({k_2}-{ k_1}\right)} }J_{mm}-\frac{\bar
t_1}{2}J_{m1}(k_1+k_2)\right)
 -(k_2-k_1)J_{m1}(J_{mm}\beta-\frac{1}{2})
 \\
 &=&
  \frac{\pl^2 f }{\pl
\bar t_1\pl \bar s_1} -J_{m1}k_1-\beta
J_{mm}J_{m1}(k_2-k_1)\eean
and 
 \bean
{\cal A}_2{\cal C}_1f
 &=&
 {\cal A}_2\left( \frac{\pl f}{\pl \bar t_1}
 -\alpha (k_2-k_1)J_{m1}\right)\\
 &=& \frac{\pl  }{\pl \bar t_1}{\cal
A}_2 f
 -{\cal A}_2(\alpha (k_2-k_1)J_{m1})\\
 &=&
 \frac{\pl ^2 }{\pl \bar t_1\pl \bar s_2}f-(k_2-k_1)\al
 J_{m1}J_{mm}
 \eean

 \bean
{\cal C}_2{\cal A}_1^+f&=&\frac{\pl   }{\pl \bar s_1
}{\cal C}_2f + {\cal C}_2 (\frac{\alpha}{2}{{
\left({k_2}-{ k_1}\right)} }J_{mm} ) \\
  &=&
    \frac{\pl^2 f }{\pl \bar s_1\pl
\bar t_2} -  {\alpha} {{ \left({k_2}-{ k_1}\right)}
}J_{m1}^2.\eean
Finally, set $\beta=0$ in each of these expressions.
That means that we can replace all ${\cal A}_i,~{\cal
C}_i$ by $\bar{\cal A}_i,~\bar{\cal C}_i$, except for
$\bar {\cal A}_2$ and
 $\bar {\cal B}_2$
, which gives an extra $\pl/\pl \beta$ and the
composition
\bean \AR_2^+\CR_1f\big\vert_{\beta
=0}&=&\left(\bar{\AR}_2-\frac{1}{2}\frac{\pl}{\pl\beta}\right)
\left(\bar{\CR}_1 -2J_{1m}\beta \frac{\pl
 }{\pl\al}\right)f\Big\vert_{\beta
=0}\\
&=&\left(\bar{\AR}_2\bar{\CR}_1+J_{1m}
\frac{\pl}{\pl\al}\right)f-\frac{1}{2}\bar{\CR}_1\frac{\pl
f}{\pl\beta}\Bigr|_{\beta=0}
 .\eean
This establishes Proposition 4.3.\qed

\remark The five operators
$\bar{\AR}_1^{\pm},~\bar{\AR}_2,~\bar{\CR}_1,
~\bar{\CR}_2$ form a Lie algebra (upon using Lemma 3.1):
$$\left[\bar{\AR}_1^{\pm} , \bar{\CR}_1\right]=0,~~
\left[\bar{\AR}_1^{+} , \bar{\AR}_1^{-}\right]=0, ~~
\left[\bar{\AR}_2  , \bar{\CR}_2 \right]=0
$$

$$\left[\bar{\AR}_1^{\pm} , \bar{\AR}_2
\right]=-(\frac{1}{2} +J_{mm})\bar{\AR}_1^{\pm}
,~~~\left[\bar{\AR}_1^{\pm} , \bar{\CR}_2
\right]=-J_{1m} \bar{\CR}_1
$$

$$\left[\bar{\AR}_2  , \bar{\CR}_1
\right]=-J_{1m} (\bar{\AR}_1^+ + \bar{\AR}_1^-)
,~~~~\left[\bar{\CR}_2  , \bar{\CR}_1
\right]=-(1+2J_{11}) \bar{\CR}_1
$$

\newpage

\section{Integrable deformations and 3-component KP}

    The following integral
\bea \lefteqn{\tau_{k_1k_2}(\bar t,\bar s,\bar
u;\alpha,\beta,\gamma^{(2)},\ldots, \gamma^{(m-1)},c
^{(1)},\dots,c ^{(m-1)},E_1\times
\ldots\times E_m)~~~~~~~~~~~~~~~~~~~~~~~~~}\no\\
 \no\\
   &&\hspace*{3cm}=\det\left(  \begin{array}{l}
       (\mu_{ij}^+)_{1\leq i \leq k_1,~1\leq j\leq
       k_1+k_2 }
       \\  \\
       (\mu_{ij}^-)_{1\leq i \leq k_2,~1\leq j\leq
       k_1+k_2 }
       \\
       \end{array}
       \right) \label{5.B}\eea
where
  \bea
   \mu^\pm_{ij}( \bar t,\bar s,\bar u,c,\gamma)
   &=&
   \int_{\prod_1^m E_i}
    (x^{(1)})^{j-1} (x^{(m)})^{i-1}
   e^{\sum_1^{\iy}\left(\bar t_kx^{(1)k}
    -\left({{\bar s_k}\atop{\bar u_k}}\right)
    x^{(m)k}\right)}
 e^{\pm \al x_m\pm \beta x_m^2}
   \no\\&&\hspace*{5cm} F_m(x^{(1)},\ldots,x^{(m)})
      \prod_{\ell=1}^m dx^{(\ell)}
  \no\\\no\\
  &=&
   \inn
 {x^{i-1}e^{ -\sum_1^{\iy}
 \left({{\bar s_k}\atop{\bar u_k}}\right)x^k}
 e^{\pm\al x \pm \beta x^2}}{x^{j-1}e^{ \sum_1^{\iy}
   \bar t_kx^k} }_m
  \label{4.C}\eea
  with regard to the inner-product
  ($m\geq 2 $)
   $$ \inn{f}{g}_m=
   \int_{\prod _1^m E_i} f(x_m)g(x_1)F_m(x_1,\ldots,x_m)dx_1
    \ldots dx_m
    ,$$
for
$$
 F_m(x_1,\ldots,x_m):=\left(\prod_1^m e^{-\frac{x_{\ell}^2}{2}}
 \right)e^{
    \sum_{ p,q\geq 1}\sum_{\ell=1}^{m-1}
    c^{(\ell)}_{pq} x_{ \ell } ^p
  x_{ \ell+1 } ^q
  +
  \sum_{{\ell}=2}^{m-1}\sum_{r=1}^{\iy}\gamma^{(\ell)}_r
      x_{ {\ell} } ^r}
    . $$

    From (\cite{AvM3}), it follows that the function
    $\tau_{k_1k_2}$ above, which is expressed as the determinant
    of a moment matrix with regard to two different weights,
    satisfies the 3-component KP
    and thus it satisfies in particular, the following
    PDE's:
\bea \frac{\pl^2}{\pl \bar s_1\pl \bar t_1}\log
\tau_{k_1k_2} &=&-
  \frac{\tau_{k_1+1,k_2}\tau_{k_1-1,k_2}}
        {\tau_{k_1,k_2}^2}
\no\\
  \frac{\pl^2}{\pl \bar u_1\pl \bar t_1}\log \tau_{k_1k_2} &=&-
  \frac{\tau_{k_1 ,k_2+1}\tau_{k_1 ,k_2-1}}
  {\tau_{k_1,k_2}^2}
\no\\
 \frac{\pl^2}{\pl \bar s_1\pl \bar u_1}\log \tau_{k_1k_2} &=&
  \frac{\tau_{k_1+1,k_2-1}\tau_{k_1-1,k_2+1}}
  {\tau_{k_1,k_2}^2}\label{equations1}
\eea
and
\bea \frac{\pl^2}{\pl \bar s_1\pl \bar t_2}\log
\tau_{k_1k_2} &=&- \frac{1}{\tau_{k_1k_2}^2} \left[
\left(\frac{\pl }{ \pl \bar t_1}
 \tau_{k_1+1,k_2}\right)
  \tau_{k_1-1,k_2}
  -
  \tau_{k_1+1,k_2}
   \left(\frac{\pl }{ \pl \bar t_1}
 \tau_{k_1-1,k_2}\right)
  \right]
 \no\\
\frac{\pl^2}{\pl \bar s_2\pl \bar t_1}\log \tau_{k_1k_2}
&=&-
 \frac{1}{\tau_{k_1k_2}^2} \left[
 \left(\frac{\pl }{ \pl \bar s_1}
 \tau_{k_1-1,k_2}\right)
  \tau_{k_1+1,k_2}
  -
  \tau_{k_1-1,k_2}
   \left(\frac{\pl }{ \pl \bar s_1}
 \tau_{k_1+1,k_2}\right)\right]
  \no\\
\frac{\pl^2}{\pl \bar s_2\pl \bar u_1}\log \tau_{k_1k_2}
&=&
 \frac{1}{\tau_{k_1k_2}^2} \left[
 \left(\frac{\pl }{ \pl \bar s_1}
 \tau_{k_1\!-\!1,k_2\!+\!1}\right)
  \tau_{k_1\!+\!1,k_2\!-\!1}
  -
  \tau_{k_1\!-\!1,k_2\!+\!1}
   \left(\frac{\pl }{ \pl \bar s_1}
 \tau_{k_1+1,k_2-1}\right)\right].\no\\
 \label{equations2}\eea
 Thus, upon taking the ratio of the first, second and third equations
 of (\ref{equations2}) and (\ref{equations1}), one finds
 $$\frac{\pl}{\pl \bar t_1}\log \frac{\tau_{k_1+1,k_2}}
                           {\tau_{k_1-1,k_2}}
     =
      \frac
     {\frac{\pl^2}{\pl \bar t_2\pl \bar s_1}\log \tau_{k_1,k_2} }
     {\frac{\pl^2}{\pl \bar t_1\pl \bar s_1}\log \tau_{k_1,k_2}}
    $$
    $$
    \frac{\pl}{\pl \bar s_1}\log \frac{\tau_{k_1+1,k_2}}
                          {\tau_{k_1-1,k_2}}
     =
     - \frac
     {\frac{\pl^2}{\pl \bar t_1\pl \bar s_2}\log \tau_{k_1,k_2}}
     {\frac{\pl^2}{\pl \bar t_1\pl \bar s_1}\log \tau_{k_1,k_2}}
   .$$
   \be
    \frac{\pl}{\pl \bar s_1}\log \frac{\tau_{k_1+1,k_2-1}}
                          {\tau_{k_1-1,k_2+1}}
     =
     - \frac
     {\frac{\pl^2}{\pl \bar u_1\pl \bar s_2}\log \tau_{k_1,k_2}}
     {\frac{\pl^2}{\pl \bar u_1\pl \bar s_1}\log \tau_{k_1,k_2}}
  .\label{equations3}\ee


\section{A PDE for the Gaussian matrices coupled in a chain, with external potential}


The purpose of this section is to prove Theorem 0.2.
Notice the probability \bea
\lefteqn{\BP_n\left(\alpha;c_1,\ldots,c_{m-1};E_1\times
\ldots\times E_m\right)
}\no\\
& =& \frac{1}{Z_n}\int_{\prod_1^m E_i^{k_1+k_2}}
\Dt_{k_1+k_2}(x^{(1)})
   \prod_{i=1}^{m}dy_i^{(\ell)}
 \no\\
   && \Dt_{k_1}(x^{(m)'})
\prod_{i=1}^{k_1}e^{-\frac{1}{2}\sum_{\ell=1}^{m}
y_i^{(\ell)2}
  + \sum_{\ell=1}^{m-1}c_{\ell}y_i^{(\ell)}y_i^{(\ell+1)}
   +\alpha y_i^{(m)}  }\no\\
 &&   \Dt_{k_2}(x^{(m)''})
  \!\!\prod_{i=k_1+1}^{k_1+k_2}\!\!
   e^{-\frac{1}{2}\sum_{\ell=1}^{m}
y_i^{(\ell)2}
  + \sum_{\ell=1}^{m-1}c_{\ell}y_i^{(\ell)}y_i^{(\ell+1)}
   -\alpha y_i^{(m)}},
  \label{5.A} \eea
  is invariant under the involution
$$\iota:~~
 \alpha \longleftrightarrow -\alpha,~~
 \beta \longleftrightarrow -\beta,~~
 u_i\longleftrightarrow s_i,~~
  k_1 \longleftrightarrow k_2
  .$$


\medskip\noindent{\it Proof of Theorem 0.1:\/} The first equality in each of the expressions
(\ref{6.3}) and (\ref{6.4}) below follows from the
expressions for $\bar {\cal A}_1^+f $ and $\bar {\cal
C}_1f $ in (\ref{Vir-identities}), whereas the second
equalities follow from (\ref{equations3}) and the third
equality from the expressions for ${\cal A}_1^+{\cal
C}_1f$, ${\cal A}_2^+{\cal C}_1f$, ${\cal C}_2{\cal
A}_1^+f$ in (\ref{Vir-identities}):
%
%
%
%
%
%
%
 \bea
 \bar{\cal C}_1\log   \frac{\tau_{k_1+1,k_2}}
                          {\tau_{k_1-1,k_2}}&=&
 \frac{\pl}{\pl \bar t_1}\log \frac{\tau_{k_1+1,k_2}}
                           {\tau_{k_1-1,k_2}}
     + {2\alpha J_{1m}} \no\\
     &=&
      \frac
     {\frac{\pl^2}{\pl \bar t_2\pl \bar s_1}\log \tau_{k_1,k_2} }
     {\frac{\pl^2}{\pl \bar t_1\pl \bar s_1}\log \tau_{k_1,k_2}}
   + {2\alpha J_{1m}}\no\\
    &=&
      \frac
     {\left(\frac{\pl^2}{\pl \bar t_2\pl \bar s_1}
     + {2\alpha J_{1m}}\frac{\pl^2}{\pl \bar t_1\pl \bar s_1}\right)\log \tau_{k_1,k_2} }
     {\frac{\pl^2}{\pl \bar t_1\pl \bar s_1}\log
     \tau_{k_1,k_2}}= -\frac{H_2^+}{F^+} \label{6.3}
   \eea
   \bea
 \bar{\cal A}_1^+\log   \frac{\tau_{k_1+1,k_2}}
                          {\tau_{k_1-1,k_2}}&=&
 \frac{\pl}{\pl \bar s_1}\log \frac{\tau_{k_1+1,k_2}}
                          {\tau_{k_1-1,k_2}}
     -\alpha J_{mm}\no\\
     &=&
     - \frac
     {\frac{\pl^2}{\pl \bar t_1\pl \bar s_2}\log \tau_{k_1,k_2}}
     {\frac{\pl^2}{\pl \bar t_1\pl \bar s_1}\log \tau_{k_1,k_2}}
   -\alpha J_{mm}
  \no\\
   &=&-\frac{H_1^+-X
   }
   {F^+}-\alpha J_{mm}  \label{6.4}
   \eea
   where $F^+,~X,~H_i^+$ can are functions of $\log \tau_{k_1k_2}$,
   which can also be expressed in terms of the actual
   probability $\BP_n$, taking into account Lemma 8.1,
\bean F^+&=& \frac{\pl^2}{\pl \bar t_1\pl \bar s_1}\log
     \tau_{k_1,k_2}=\bar{\cal A}_1^+\bar{\cal C}_1 \log
\tau_{k_1,k_2}+k_1J_{1m}=
 \bar{\cal A}_1^+\bar{\cal C}_1 \log \BP_{n}+k_1J_{1m}
\\
 X &=&
 \frac{1}{2}\bar{\cal C}_{ 1}\frac{\pl  }{\pl \beta }
    \log \tau_{k_1k_2}
   =\frac{1}{2}\bar{\cal C}_{ 1}\frac{\pl  }{\pl \beta }
    \log \BP_n
%
\\
 H_1^{+}&=&
  \frac{\pl^2}{\pl \bar t_1\pl \bar s_2}\log \tau_{k_1,k_2}
  +
   \frac{1}{2}\frac{\pl}{\pl \beta}\bar{\cal C}_1\log \tau_{k_1,k_2}\\
  &=&  \left( \bar{\cal
A}_2\bar{\cal C}_1 + J_{1m}\frac{\pl}{\pl \alpha}\right)
  \log \tau_{k_1,k_2} 
-\alpha ( k_1-k_2)J_{1m}J_{mm}
    \\
     &=& \left( \bar{\cal
A}_2\bar{\cal C}_1+ J_{1m}\frac{\pl}{\pl \alpha}\right)
  \log \BP_{n} 
- 2\alpha k_1 J_{1m}J_{mm}+\frac{k_1k_2}{\alpha}J_{1m}
\eean
\bean
 H_2^{+}&=&
  - \left(\frac{\pl^2}{\pl \bar t_2\pl \bar s_1}
     + {2\alpha J_{1m}}\frac{\pl^2}{\pl \bar t_1\pl \bar s_1}\right)
     \log \tau_{k_1,k_2} \\
     &=&
   - \left( \bar{\cal
C}_2 + 2\alpha J_{1m}\bar{\cal C}_1\right)\bar{\cal
A}_1^+ \log \tau_{k_1,k_2}
  -\alpha (k_1+k_2)(J_{1m})^2
 \\
  &=&
   - \left( \bar{\cal
C}_2 + 2\alpha J_{1m}\bar{\cal C}_1\right)\bar{\cal
A}_1^+ \log \BP_{n} \eean
   Thus, $\bar {\cal A}_1^+$ of the right hand side of (\ref{6.3}) must equal $\bar {\cal
   C}_1$ of the right hand side of (\ref{6.4}) and furthermore noticing that $\bar {\cal
   C}_1\al J_{mm}=0$, one finds:
$$
\bar {\cal A}_1^+\frac{H_2^+}{F^+}=\bar {\cal
C}_1\left(\frac{H_1^+}{F^+}-\frac{X}{F^+}\right)
 $$
or in Wronskian notation
  \be  \left\{X,
F^+\right\}_{\bar{\CR}_1}=\{H^+_1,F^+\}_{\bar{\CR}_1}-\{H_2^+,F^+\}_{\bar{\AR}_1}
=G^+
 \label{6.5}\ee
 and
 \be -\left\{X,   F^- \right\}_{\bar{\cal C}_{ 1}} =   \left\{
 H_1^-,F^- \right\}_{\bar{\cal C}_{ 1}}
      -\left\{H^-_2,   F^- \right\}_{\bar{\cal B}_{ 1}}
    =G^- \label{6.6}\ee
     upon applying the involution $\iota$ to the first
     equation.
     Equations (\ref{6.5}) and (\ref{6.6}) yield a linear
  system of equations in
   $
   X~~\mbox{and}
   ~~ \bar{\cal C}_{ 1}X,
   $
leading to
  \bean
    X&=&
    \frac{G^- F^+ + G^+F^-}
    {-F^-(\bar{\cal C}_{ 1}F^+)+F^+(\bar{\cal C}_{ 1}F^-)}
\\
 \bar{\cal C}_{ 1}X &=&
    \frac{G^- (\bar{\cal C}_{ 1}F^+) + G^+(\bar{\cal C}_{ 1}F^-)}
    {-F^-(\bar{\cal C}_{ 1}F^+)+F^+(\bar{\cal C}_{ 1}F^-)}
\eean
  Subtracting the second equation from
$\bar{\cal C}_{ 1}$ of the first equation yields the
following:
  \bean&&
\Bigl(F^+\bar{\cal C}_{ 1}G^-+F^-\bar{\cal C}_{ 1}G^+
\Bigr)
\Bigl(F^+\bar{\cal C}_{ 1}F^- -F^-\bar{\cal C}_{ 1}F^+ \Bigr) \\
  &&-
  \Bigl(F^+ G^- +F^- G^+
\Bigr) \Bigl(F^+\bar{\cal C}_{ 1}^2F^- -F^-\bar{\cal
C}_{ 1}^2F^+ \Bigr) =0.
 \eean

The second way of expressing these equations is to write
a system of 4 equations, consisting of the system
(\ref{6.5}) and (\ref{6.6}) and that same system acted
upon by $\bar{\cal C}_1$:
\bean
0&=&-\left\{X,F^+\right\}_{\bar{\cal C}_1}~~+~G^+ \\
0&=&-\left\{X,F^-\right\}_{\bar{\cal C}_1}~~-~G^- \\
0&=&-\bar{\cal C}_1\left\{X,F^+\right\}_{\bar{\cal C}_1}+\bar{\cal C}_1G^+ \\
0&=&-\bar{\cal C}_1\left\{X,F^-\right\} _{\bar{\cal
C}_1}-\bar{\cal C}_1G^- ,\eean or, in matrix notation,
$$
\left( \begin{array}{cccc}
   G^+& \bar{\cal C}_1F^+  & - F^+ &0\\
   -G^-& \bar{\cal C}_1F^-  & - F^- &0\\
   \bar{\cal C}_1 G^+ &\bar{\cal C}_1^2F^+&0&-F^+\\
    -\bar{\cal C}_1 G^- &\bar{\cal C}_1^2F^-&0&-F^-\\
  \end{array}\right)
   \left(
      \begin{array}{c}
      1\\
      X\\
      \bar{\cal C}_1 X\\
      \bar{\cal C}_1^2X
      \end{array}
\right)=0
$$
and thus the matrix must be singular, leading to the
second formulation and ending the proof of Theorem
0.1.\qed



    \section{The PDE for the transition probability of the Pearcey process}

Given $E_{\ell}:=
\bigcup^r_{i=1}[x^{(\ell)}_{2i-1},x^{(\ell)}_{2i}]\subset
\BR$, define the ``{\em space}" and ``{\em time}"
operators $ {\cal X}_k$ and $ {\cal T}_k$ related to
Brownian motion, together with a mixed ``{\em
space-time}" operator $\tilde{\cal X}_{-1}$,
 $$
 \hspace*{-1cm} {\cal
X}_k:=\sum_{\ell}
\sum_{i=1}^{2r_{\ell}}(x_i^{(\ell)})^{k+1}\frac{\pl}{\pl
x_i^{(\ell)}},~~~{\cal
T}_j=\sum_{\ell}s_{\ell}^{j+1}\frac{\pl}{\pl
s_{\ell}},~~~ \tilde{\cal
X}_{-1}=\sum_{\ell}s_{\ell}\sum_{i=1}^{2}\frac{\pl}{\pl
x_i^{(\ell)}}.$$ We shall also need the intermediate
operators
 $$
 {\cal
X}_k(x^{(\ell)}):=
\sum_{i=1}^{2r_{\ell}}(x_i^{(\ell)})^{k+1}\frac{\pl}{\pl
x_i^{(\ell)}} .$$
Define a new function $\BQ_z$
\be
\log\BP_n(\al,c_1,\ldots,c_{n-1};b^{(1)},\ldots,b^{(m)})
\Bigr|_{n=\frac{2}{z^4}}
=\BQ_z(s_1,\ldots,s_m;x^{(1)},\ldots,x^{(m)})
\label{7.2}\ee
by means of the change of variables, defined earlier
\bean
  \al&=&\frac{2 }{z}\frac{\sqrt{s_m-s_{m-1}}}{\sqrt{(1-s_mz^2)(1-s
_{m-1}z^2)}}
 \\
c_j &=&\left\{\begin{array}{ll}
\sqrt{\frac{(1+s_1z^2)(s_3-s_2)}{(1+s_2z^2)(s_3-s_1)}}&\mbox{for~}j=1\\
\\
\sqrt{\frac{(s_j-s_{j-1})(s_{j+2}-s_{j+1})}{(s_{j+1}-s_{j-1})
(s_{j+2}-s_j)}}&\mbox{for~}2\leq j\leq
m-2\\
\\
\sqrt{\frac{(s_{m-1}-s_{m-2})(1-s_mz^2)}{(s_m-s_{m-2})(1-s_{m-1}z^2)}}&\mbox{for~}j=m
-1
\end{array}\right.\\
\\
 b_i^{(\ell)}&=& \left\{\begin{array}{l}
2x_i^{(1)}\sqrt{\frac{1+s_2z^2}{(1+s_1z^2)(s_2-s_1)}}
 ,~\mbox{~~for $\ell=1$}\\
 2x_i^{(\ell)}
 \sqrt{\frac{s_{\ell+1}-s_{\ell-1}}{(s_{\ell}-s_{\ell-1})
 (s_{\ell+1}-s_{\ell})}}
  ,~\mbox{~~for $2\leq \ell\leq m-1$} \\
2x_i^{(m)}\sqrt{\frac{1-s_{m-1}z^2}{(s_m-s_{m-1})(1-s_mz^2)}}
 ,~\mbox{~~for $\ell=m$}
\end{array}\right.
\eean \be \label{7.3}\ee
  with inverse map
  \bea
s_i&=&s_i(\al,c;z)=\frac{1}{z^2}\frac{\al^2z^4\frac{J_{m1}J_{mi}}{J_{1i}}
+2 }{\al^2z^4\frac{J_{m1}J_{mi}}{J_{1i}}-2 } \no\\\no \\
 x^{(i)}_k&=&b^{(i)}_k U_i(s(\al,c;z);z)
,\label{7.4}
 \eea
 with
 $$
  U_i(s(\al,c;z);z)
  =
  \frac{-1}{J_{mi}}\left(\frac{\al z~\frac{J_{m1}J_{mi}}{J_{1i}}}
 {\al^2z^4\frac{J_{m1}J_{mi}}{J_{1i}}-2}\right)
 .$$
Then
\bean \frac{\pl x_k^{(i)}}{\pl\al}&=&\frac{\pl
(b_k^{(i)}U_i)}{\pl\al}=b_k^{(i)}\frac{\pl
U_i}{\pl\al}=b_k^{(i)}U_i\frac{\pl \log
U_i}{\pl\al}=x_k^{(i)}\frac{\pl \log U_i}{\pl\al}\\
\\
\frac{\pl x_k^{(i)}}{\pl c_j}&=&\frac{\pl
(b_k^{(i)}U_i)}{\pl c_j}=b_k^{(i)}\frac{\pl U_i}{\pl
c_j}=b_k^{(i)}U_i\frac{\pl \log U_i}{\pl
c_j}=x_k^{(i)}\frac{\pl \log U_i}{\pl c_j} .\eean
Thus, setting (\ref{7.4}) in the right hand side of
(\ref{7.2}) and using the chain rule, one computes
 \bean
\frac{\pl\log\BP_n}{\pl\al}\Bigr|_{n=\frac{2}{z^4}}
 &=&\sum^m_{i=1}\frac{\pl
s_i}{\pl\al}\frac{\pl \BQ_z}{\pl s_i}+\sum^m_{i=1}
\sum^{2r_i}_{k=1}\frac{\pl x_k^{(i)}}{\pl\al}
 \frac{\pl \BQ_z}{\pl x_k^{(i)}}\\
\\
&=&\left(\sum^m_{i=1}\frac{\pl s_i}{\pl\al}\frac{\pl
}{\pl s_i}+\sum^m_{i=1}\frac{\pl \log
U_i}{\pl\al}\sum^{2r_i}_{k=1}x_k^{(i)}\frac{\pl }{\pl
x_k^{(i)}}\right)\BQ_z\\
\\
&=&\sum^m_{i=1}\left(\frac{\pl s_i}{\pl\al} \frac{\pl
}{\pl s_i}+\frac{\pl \log
U_i}{\pl\al}\XR_0(x^{(i)})\right)\BQ_z \eean

\bean
 \frac{\pl\log\BP_n}{\pl
c_j}\Bigr|_{n=\frac{2}{z^4}}
 &=&\sum^m_{i=1}\frac{\pl
s_i}{\pl c_j}\frac{\pl\BQ_z}{\pl s_i}+\sum^m_{i=1}
\sum^{2r_i}_{k=1}\frac{\pl x_k^{(i)}}{\pl c_j}\frac{\pl \BQ_z}{\pl x_k^{(i)}}\\
\\
&=&\left(\sum^m_{i=1}\frac{\pl s_i}{\pl c_j}\frac{\pl
}{\pl s_i}+\sum^m_{i=1}\frac{\pl \log U_i}{\pl
c_j}\sum^{2r}_{k=1}x_k^{(i)}\frac{\pl }{\pl
x_k^{(i)}}\right)\BQ_z\\
\\
&=&\sum^m_{i=1}\left(\frac{\pl s_i}{\pl c_j} \frac{\pl
}{\pl s_i}+\frac{\pl \log U_i}{\pl
c_j}\XR_0(x^{(i)})\right)\BQ_z
 \eean
 \bean
\DR_{-1}(E_{i})\log\BP_n\Bigr|_{n=\frac{2}{z^4}}
 &=&\sum^{2r_i}_{k=1}\frac{\pl \log\BP_n}{\pl b_k^{(i)}}
  \Bigr|_{n=\frac{2}{z^4}}\\
\\
&=& \sum^{2r_i}_{k=1}\frac{\pl x_k^{(i)}}{\pl
b_k^{(i)}}\frac{\pl \BQ_z}{\pl
x_k^{(i)}}\\
\\
&=&\left(U_i\sum^{2r_i}_{k=1}\frac{\pl }{\pl x_k^{(i)}}
\right)
 \BQ_z\\
\\
&=&U_i\XR_{-1}(x^{(i)})\BQ_z
 \eean
 \bean \DR_0(E_{i})\log\BP_n\Bigr|_{n=\frac{2}{z^4}}
&=&\sum^{2r_i}_{k=1}b_k^{(i)}\frac{\pl \log\BP_n}{\pl
b_k^{(i)}}\Bigr|_{n=\frac{2}{z^4}}\\
\\
&=& \sum^{2r_i}_{k=1}b_k^{(i)}\frac{\pl x_k^{(i)}}{\pl
b_k^{(i)}}\frac{\pl \BQ_z}{\pl x_k^{(i)}} \\
\\
&=&\left(\sum^{2r_i}_{k=1}x_k^{(i)}\frac{\pl }
 {\pl x_k^{(i)}} \right)\BQ_z\\
\\
&=& \XR_0(x^{(i)})\BQ_z
 .\eean
Using the information above, one computes
  ($\vr=\pm 1$)
\bea \lefteqn{\bar\AR_1^{\vr} \log\BP_n
(\al,c_1,\ldots,c_{n-1};b^{(1)},\ldots,b^{(m)})
\Bigr|_{n=\frac{2}{z^4}}}
  \no\\
&=&-\frac{1}{2}\left(\sum^m_{j=1}J_{mj}\DR_{-1}(E_{j})
+\vr\frac{\pl }{\pl\al}\right)
   \log\BP_n\Bigr|_{n=\frac{2}{z^4}}
    \no\\
  \no\\
&=&-\frac{1}{2}\sum^m_{j=1}\left(J_{mj}U_j\XR_{-1}(x^{(j)})
+\vr\frac{\pl s_j}{\pl\al} \frac{\pl }{\pl
s_j}+\vr\frac{\pl \log U_j}{\pl\al}
\XR_0(x^{(j)})\right)\BQ_z
 \no\\
\no\\
&=:&A_1^{\vr}\BQ_z
 \label{7.A1}
 \\\no\\
\bar\CR_1\log\BP_n\Bigr|_{n=\frac{2}{z^4}}
 &=&\sum^m_{j=1}J_{1j}\DR_{-1}(E{j})\log\BP_n
 \Bigr|_{n=\frac{2}{z^4}}
  \no\\
  \no\\
&=&\left(\sum^m_{j=1}J_{1j}U_j\XR_{-1}(x^{(j)})\right)\BQ_z
  \no\\
  \no\\
&=:&C_1\BQ_z
  \label{7.C1}\eea

\bea
\lefteqn{\bar\AR_2\log\BP_n\Bigr|_{n=\frac{2}{z^4}}
 }\no\\&=&
  \frac{1}{2}
 \left(
  \DR_0(b^{(m)})- \alpha\frac{\pl}{\pl \alpha}
  - c_{m-1}\frac{\pl}{\pl c_{m-1}}
  \right)\log\BP_n\Bigr|_{n=\frac{2}{z^4}}
  \no\\
  &=&-\frac{1}{2}\sum^{m}_{j=1}
\left(\begin{array}{l} \left( c_{m-1}\frac{\pl s_j}{\pl
c_{m-1}}+\al\frac{\pl s_j}{\pl\al} \right)
\frac{\pl}{\pl s_j}\no\\ \no \\
 +\left( \left(
c_{m-1}\frac{\pl\log U_j}{\pl c_{m-1}}+\al\frac{\pl\log
U_j}{\pl\al} \right) -\dt_{jm} \right) \XR_0(x^{(j)})
\end{array}
\right)
\BQ_z \no\\
\no\\
&=:&A_2\BQ_z
\label{7.A2}\eea

  \bea
\bar\CR_2\log\BP_n\Bigr|_{n=\frac{2}{z^4}}&=&
  \left(-
  \DR_0(b^{(1)}) + c_1\frac{\pl}{\pl c_1}
 \right)\log\BP_n\Bigr|_{n=\frac{2}{z^4}}
   \no\\
 &=&\sum^m_{j=1}\left(c_1\frac{\pl s_j}{\pl
c_1}\frac{\pl }{\pl s_j}+ \left(c_1\frac{\pl \log
U_j}{\pl c_1}-\dt_{j1}\right)\XR_0(x^{(j)})\right)\BQ_z
 \no\\
  \no\\
&=:&C_2\BQ_z
\label{7.C2}\eea


\begin{lemma}
  The following expansions hold near $z \sim 0$:
  \bean
  \al \frac{\pl s_i}{\pl \al}&=&
  \frac{1}{z^2}-s_i^2z^2
  \\
 \frac{\pl s_i}{\pl \al}
  &=& \frac{1}{2z\sqrt{s_m-s_{m-1}}}\left(1
  -\frac{z^2}{2}(s_{m-1}+s_m)+O(z^4)\right)
  \eean
  \bean
  \al\frac{\pl}{\pl \al}\log U_i (s(\al,c;z),z)&=&-z^2s_i
  \\
   \frac{\pl}{\pl \al}\log U_i (s(\al,c;z),z)&=&O(z^3)
  \eean
  \bean
  c_{ 1}\frac{\pl }{\pl c_{ 1}}
  s_i&=&\frac{1}{2(s_2-s_{ 1})}\left(\frac{1}{z^4}+\frac{s_1+s_{2}
  -2s_i}{z^2}+O(1) \right)\\
  c_{m-1}\frac{\pl }{\pl c_{m-1}}
  s_i&=&
   \frac{1}{s_m-s_{m-1}}\left(\frac{1}{z^4}+\frac{s_m+s_{m-1}
  -2s_i}{z^2}+O(1)\right)
  \eean
  \bean
  c_{ 1}\frac{\pl }{\pl c_{ 1}}
  \log U_i&=&-\frac{1}{2(s_2-s_{ 1})z^2}+O(1)
  \\
  c_{ m-1}\frac{\pl }{\pl c_{m-1}}
  \log U_i&=&-\frac{1}{2(s_m-s_{m-1})z^2}+O(1)
  \eean
  \bean
  J_{mi} U_i&=&\frac{1}{4\sqrt{s_m-s_{m-1}}}
  \left(-\frac{1}{z^2}+\frac{1}{2} (s_m+s_{m-1}-2s_i)+O(z^2)\right)
  \\
  J_{1i} U_i&=&\frac{1}{4\sqrt{s_2-s_{1}}}
  \left(-\frac{1}{z^2}-\frac{1}{2} (s_1+s_{2}-2s_i)\right.
  \\&&~~~~~~~~~~~~~~~~~
  \left.+\frac{z^2}{8}
  ((s_2-s_1)^2
  +4s_1(s_1+s_2))+O(z^4)\right)
  \eean
  $$
   \al=\frac{2 }{z}\sqrt{s_m-s_{m-1}}~
   (1+\frac{z^2}{2}(s_m+s_{m-1})+O(z^4))
 . $$
\end{lemma}

\proof These series follow from expanding the
expressions given in Lemma 3.3 near $z\sim 0$. \qed

\begin{lemma} 
The operators $A^{\vr}_1, ~A_2, ~, C_1, ~C_2$, as
defined in (\ref{7.A1}), (\ref{7.A2}), (\ref{7.C1}),
(\ref{7.C2}) and as acting on the function
  $\BQ_z(s_1,\ldots,s_m;x^{(1)},\ldots,x^{(m)})$ admit the
following expansions in $z\thicksim 0$, in terms of the
operators ${\cal X}_k,~{\cal T}_k, ~\bar{\cal X}_{-1}$,
\bean
  A^{\vr}_1&=&\frac{1}{8\sqrt{s_m-s_{m-1}}}
\left\{\begin{array}{l}
\frac{1}{z^2}\XR_{-1}-\frac{2\vr}{z}{\BT}_{-1}+
 \left(\tilde\XR_{-1}-\frac{(s_{m-1}+s_m)}{2}\XR_{-1}\right)\\
+\vr z(s_{m-1}+s_m){\cal T}_{-1}+{O}(z^2)
\end{array}\right\}
 \\
 \\
C_1&=&\frac{1}{4\sqrt{s_2-s_1}}\left\{\begin{array}{l}
-\frac{1}{z^2}{\cal X}_{-1}+\left(\tilde\XR_{-1}-
\frac{(s_1+s_2)}{2}\XR_{-1}\right)
\\
+z^2\left(\frac{(s_2-s_1)^2}{8}
\XR_{-1}+\frac{(s_1+s_2)}{2}\tilde\XR_{-1}\right)+{ O}(z^4)
\end{array}\right\}
 \\
 \\
 A_2&=&\frac{1}{4(s_m-s_{m-1})}\left\{ -\frac{{\cal
T}_{-1}}{z^4}+\frac{1}{z^2}(\XR_0+2{\cal T}_0+(s_{m-1}-3
s_m){\cal T}_{-1})+ {  O}(1) \right\}
 \\
 \\
 C _2&=&\frac{1}{2(s_2-s_1)}\left\{ \frac{{\cal
T}_{-1}}{z^4}-\frac{1}{z^2}(\XR_0+2{\cal
T}_0-(s_1+s_2){\cal T}_{-1})+ {  O}(1) \right\} .\eean
\end{lemma}

\proof These formulae are an immediate consequence of
formulae (\ref{7.A1}), (\ref{7.A2}), (\ref{7.C1}),
(\ref{7.C2}), combined with the series in Lemma 7.1.\qed

\bigbreak

%

  \begin{lemma}  In the new coordinates $s_1,\ldots,s_m,x^{(1)},\ldots,x^{(m)}$,
   the quantities appearing in the basic equation (\ref{BasicEquation}) for $\BP_n$
   admit the following expansions in $z\sim 0$, setting $k=n/2=1/z^4$:
     \newcommand{\ER}{{\cal E}}
\bea \lefteqn{F^{\vr}=
 { A}_1^\vr{  C}_1  \BQ_{z}+\frac{J_{1m}}{z^4}}
\no\\ 
&&=\frac{1}{32\sqrt{(s_2-s_1)(s_m-s_{m-1})}} \no\\
&&\left\{
\begin{array}{l}
-\frac{16}{z^6}-\frac{1}{z^4}(\XR^2_{-1}\BQ_z
 +8(s_1+s_2-s_{m-1}-s_m))\\ \\
+~\frac{2\vr}{z^3}\XR_{-1}{\cal T}_{-1}\BQ_z
+\frac{1}{z^2}\left[
\frac{1}{2}(s_{m-1}+s_m-s_1-s_2)\XR^2_{-1}\BQ_z\right.\\
\\
\left.+~2((s_2-s_1)^2+(s_m-s_{m-1})^2+2(s_1+s_2)(s_{m-1}+s_m)
 )\right]\\
\\
-\frac{\vr}{z}((s_{m-1}+s_m-s_1-s_2)\XR_{-1}
+2\tilde\XR_{-1}){\cal T}_{-1}\BQ_z+{  O}(1)
\end{array}\right\}
\eea

\bea
\lefteqn{\bar\CR_1F^{\vr}=\frac{1}{128(s_2-s_1)\sqrt{s_m-s_{m-1}}}
}\no\\
&&\left\{\begin{array}{l}
\frac{1}{z^6}\XR^3_{-1}\BQ_z-\frac{2\vr}{z^5}\XR^2_{-1}
{\cal T}_{-1}\BQ_z\\
+ \frac{1}{z^4}\left((s_1 +s_2-\frac{(s_{m-1}+s_m)}{2})
\XR_{-1} 
-\tilde\XR_{-1}\right)\XR^2_{-1}\BQ_z\\
+\frac{\vr}{z^3}( (s_{m-1}+s_m-2(s_1+s_2))\XR_{-1}
+~4\tilde\XR_{-1}) \XR_{-1}{\cal T}_{-1}\BQ_z  +{
O}\left(\frac{1}{z^2}\right)
\end{array}\right\}\no\\
\eea

\bea  \lefteqn{\bar\CR^2_1F^{\vr}=
\frac{1}{512(s_2-s_1)^{3/2}\sqrt{(s_m-s_{m-1})}}
 } \\
 &&~~~~~~~~~~~~
\left\{-\frac{\XR^4_{-1}\BQ_z}{ z^8}+
\frac{2\vr}{z^7}\XR^3_{-1}{\cal T}_{-1}\BQ_z+{
O}\left(\frac{1}{z^6}\right)\right\}
 \nonumber\\
\lefteqn{\bar\AR_1F^{\vr}=\frac{1}{256(s_m-s_{m-1})\sqrt{(s_2-s_{1})}}
}\no\\
&&~~~~~~~~~~~ \left\{-\frac{\XR^3_{-1}\BQ_z}{z^6} +
\frac{4\vr}{z^5}\XR^2_{-1}{\cal T}_{-1}\BQ_z+{
O}\left(\frac{1}{z^4}\right)\right\} \eea

\bea
 \lefteqn{H_1^{\vr}=\left(  {A}_2{C}_1+\vr J_{1m}\frac{\pl}{\pl
\alpha}\right)
   \BQ_{z} 
- \vr\frac{J_{1m}}{z^4} \left(2\alpha
 J_{mm}-\frac{ 1}{z^4\alpha} \right)
 }
 \\
 &=&\frac{1}{16(s_m-s_{m-1})\sqrt{s_2-s_1}}
 \no\\
 &&
 \left\{\begin{array}{l}
-\frac{20\vr}{z^9}-\frac{2\vr}{z^7}(5(s_1+s_2)-10s_{m-1}+6s_m)
%
+\frac{1}{z^6}\XR_{-1}{\cal T}_{-1}\BQ_z
\no\\
+\frac{\vr}{2z^5}\left(5(s_2-s_1)^2+
 24s_{m-1}s_m+4(s_1+s_2)(5s_{m-1}-3s_m)\right)\\
+\frac{1}{z^4}\left(\begin{array}{l}
 2\XR_{-1}-\tilde\XR_{-1}{\cal T}_{-1}
 -2\XR_{-1}{\cal T}_0-\XR_{-1}\XR_0\\
  +
\left(\frac{(s_1+s_2)}{2}-s_{m-1}+3s_m\right)
\XR_{-1}{\cal T}_{-1}\end{array}\right)\BQ_z
+{O}\left(\frac{1}{z^3}\right)
\end{array}\right\}\no\\
\eea

\bea
 \lefteqn{\bar\CR_1H_1^{\vr}=\frac{1}{64(s_2-s_1)(s_m-s_{m-1})}}
\no\\
 &&
\left\{\begin{array}{l} -\frac{1}{z^8}\XR_{-1}^2{\cal
T}_{-1}\BQ_z
+\frac{1}{z^6}\left(\begin{array}{l}
2\tilde\XR_{-1}\XR_{-1}{\cal T}_{-1} +2\XR^2_{-1}{\cal
T}_0
+\XR_0\XR^2_{-1} \\
+~(s_{m-1}-3s_m-s_1-s_2)\XR_{-1}^2{\cal T}_{-1}
\end{array}\right)\BQ_z\\
+~{O}\left(\frac{1}{z^5}\right)
\end{array}\right\}
\no\\
\eea

\bean  -H_2^{\vr}&=&(C_2+2\vr\al J_{1m}C_1)A_1^\vr
\BQ_z \\
&=&\frac{1}{16(s_2-s_1)\sqrt{s_m-s_{m-1}}}
\\
&&\left\{\begin{array}{l}
\frac{\vr}{z^7}\XR^2_{-1}\BQ_z-\frac{1}{z^6}
 \XR_{-1}{\cal
T}_{-1}\BQ_z\\
+\frac{\vr}{z^5} (-2{\cal T}^2_{-1}+
(s_1+s_2-\frac{(s_{m-1}+s_m)}{2} )
\XR^2_{-1})\BQ_z\\
+\frac{1}{z^4} \left(\!\begin{array}{l}  \XR_{-1}
-\XR_0\XR_{-1}-2\XR_{-1}{\cal T}_0
 +3\tilde\XR_{-1}{\cal T}_{-1}\\
  +
\left(\frac{1}{2}(s_{m-1}+s_m)-s_1-s_2\right)
 \XR_{-1}{\cal T}_{-1}
 \end{array}\!\!\right){\BQ}_z\\
+{ O}\left(\frac{1}{z^3}\right)
\end{array}\!\!\right\}
\eean

\bean -\bar A_1^{\vr}H_2^{\vr}&=&\frac{1}{128
(s_2-s_1)(s_m-s_{m-1})} \\
&&\left\{\begin{array}{l}
\frac{\vr}{z^9}\XR_{-1}^3\BQ_z-\frac{3}{z^8}\XR_{-1}^2
 {\cal T}_{-1} \BQ_z
  \\
  +\frac{\vr}{z^7}(\tilde\XR_{-1}
+~(s_1+s_2-s_{m-1}-s_m)\XR _{-1})\XR_{-1}^2\BQ_z \\
+\frac{1}{z^6} (-\XR^2_{-1}\XR_0-2\XR_{-1}^2{\cal T}_0
 +2\tilde\XR_{-1}\XR_{-1}{\cal T}_{-1}\\
+~4{\cal T}^3_{-1}+3(s_m+s_{m-1}-s_1-s_2)\XR^2_{-1}
{\cal T}_{-1})\BQ_z\\
+{ O}\left(\frac{1}{z^5}\right)
\end{array}\right\}
\eean

\bea G^{\vr}&=&\left\{ H_1^\vr,F^\vr\right\}_{{ C}_{1}}
      -\left\{ H^\vr_2,   F^\vr \right\}_{{ A}_{ 1}}\\
     &=&\frac{1}{2\cdot(32)^2((s_2-s_1)(s_m-s_{m-1}))^{3/2}} \\
\nonumber\\
& &\left\{\begin{array}{l}
\frac{12\vr}{z^{15}}\XR^3_{-1}\BQ_z \\
+\frac{\vr}{z^{13}}
(-28\tilde\XR_{-1}
+(18(s_1+s_2-s_{m-1})+14s_m)\XR_{-1})\XR^2_{-1}\BQ_z \nonumber\\
+\frac{1}{z^{12}} \left(\begin{array}{l}
\frac{1}{2}\{\XR_{-1}{\cal
T}_{-1}\BQ_z,\XR_{-1}^2\BQ_z\}_{\XR_{-1}}
 \\
  -16({\cal T}_0-1+\frac{1}{2}{\cal
X}_0) {\cal X}_{-1}^2
\BQ_z
\\
 +32(\tilde\XR_{-1}\XR_{-1}
  -
  {\cal T}^2_{-1}){\cal T}_{-1}\BQ_z
\end{array}\right) \\
+{O}\left(\frac{1}{z^{11}}\right)
\end{array}\right\}
\eea

\end{lemma}
\proof This is done by straightforward computation,
using Lemma 7.2 and the series for the constants
appearing in $F^{\vr},~H_1^\vr,~H_2^\vr,~G^\vr$, namely
\bean
\lefteqn{J_{m1}(c)=-\frac{1}{2z^2\sqrt{(s_2-s_1)(s_m-s_{m-1})}}\left(1+\frac{z^2}{2}
(s_1+s_2-s_{m-1}-s_m)\right.}\\
&
&-\left.\frac{z^4}{8}((s_2-s_1)^2+(s_m-s_{m-1})^2+2(s_1+s_2)(s_{m-1}+s_m))+{\bf
O}(z^4)\right) \eean
$$
(\al
J_{m1})=-\frac{1}{\sqrt{s_2-s_1}}\left(\frac{1}{z^3}+\frac{1}{2z}(s_1+s_2)-
\frac{z}{8}(s_2-s_1)^2+{O}(z^3)\right)
$$
\bean
\lefteqn{\left(\frac{k^2J_{m1}}{\al}-2\al kJ_{1m}J_{mm}\right)}\\
&=&\frac{1}{4z^9(s_m-s_{m-1})\sqrt{s_2-s_1}}
\\ \\
&&\left(
\begin{array}{l}
-5-\frac{z^2}{2}(5(s_1+s_2)-10s_{m-1}+6s_m)\\
+\frac{z^4}{8}(5(s_2-s_1)^2+24s_{m-1}s_m+4(s_1+s_2)(5s_{m-1}-3s_m))+{
O}(z^6)
\end{array}\right)
\eean

\medskip\noindent{\it Proof of Theorem 0.2:\/}
In terms of $T_i$'s, defined by
 \bean
T_1&:=&F^+\bar\CR_1G^-+F^-\bar \CR_1G^+=\bar\CR_1T_3-T_5 \\
T_2&:=&\{F^-,F^+\}_{\bar\CR_1}\\
T_3&:=&F^+G^-_1+F^-G^+_1  \\
 T_4&:=&\bar\CR_1T_2=F^+\bar\CR^2_1F^--F^-\bar\CR^2_1F^+\\
 T_5&:=&
G^-\bar\CR_1F^++G^+\bar\CR_1F^-,
 \eean
the fundamental equation (\ref{BasicEquation}) can be
written
$$
0=-T_1\cdot T_2+T_3\cdot T_4=\{T_2,T_3\}_{\bar{\cal
C}_1}+T_2T_5
$$
 meaning that $T_2,~T_3$ and $T_5$ only are needed,
 which one checks to have the following series in $z$,
 using Lemma 7.3,
\bea  T_2
&=&-\frac{\XR^2_{-1}{\cal
T}_{-1}\BQ_z}{64(s_2-s_1)^{3/2}(s_m-s_{m -1})z^{11}}+{
O}\left(\frac{1}{z^{9}}\right)
\no\\
 T_3
&=&\frac{1}{128(s_2-s_1)^2(s_m-s_{m-1})^2z^{18}}\nonumber\\
& &\left\{\begin{array}{l} \frac{1}{16}\{\XR_{-1}{\cal
T}_{-1}\BQ_z,\XR^2_{-1}\BQ_z\}_{\XR_{-1}}
%
%
%
\no\\+({\cal T}_0-1+\frac{1}{2}\XR_0)\XR^2_{-1}\BQ_z
 \no\\
  -2(\tilde\XR_{-1}\XR_{-1}-{\cal T}^2_{-1}){\cal T}_{-1}\BQ_z
  \no\\-\frac{3}{32}(\XR_{-1}^2 \BQ_z)( \XR_{-1}^2 {\cal T}_{-1}\BQ_z) \\
 \end{array}\right\}+O\left(\frac{1}{z^{17}}\right) %
\no\\ \no\\
T_5
 &=&\frac{3(\XR^2_{-1}{\cal
T}_{-1}
\BQ_z)(\XR_{-1}^3\BQ_z)}{128^2(s_2-s_1)^{5/2}(s_m-s_{m-1})^{2}
 z^{20}}+{O}\left(\frac{1}{z^{19}}\right). \label{finalexpansion}\eea
Then remembering from Lemma 7.2,
 $$
C_1=\frac{-1}{4\sqrt{s_2-s_1}}
 \left(\frac{1}{z^2}{\cal X}_{-1}+O(1)\right)
, $$
one easily computes (letting $T_i^{(0)}$ be the leading
coefficient of $T_i$ in (\ref{finalexpansion}))
 \bean
  \lefteqn{\{T_2,T_3\}_{{  C}_1}}\no\\
   &=&
\frac{1}{z^{31}}\left\{T_2^{0},T_3^{0}\right\}_{{ C}_1}
+O\left(\frac{1}{z^{30}}\right) \\
&=&
  \frac{ 1}{32^3(s_2\!-\!s_1)^4(s_m- s_{m-1})^3z^{31}}\\
&&\Biggl\{\XR^2_{-1}{\cal T}_{-1}\BQ_z,\frac{1}{16}\{
\XR_{-1}{\cal T}_{-1} \BQ_z,\XR^2_{-1}
\BQ_z\}_{\XR_{-1}} +\left({\cal
T}_0-1+\frac{\XR_0}{2}\right)\XR_{-1}^2 \BQ_z
 \\
& & \hspace*{.2cm}
 - 2(\tilde\XR_{-1}\XR_{-1}-{\cal T}^2_{-1}){\cal
T}_{-1} \BQ_z -\frac{3}{32}
  (\XR_{-1}^2{\cal T}_{-1}\BQ_z)(\XR_{-1}^2\BQ_z)
     \Biggr\}_{\XR_{-1}}\!\!\! +\!\!
{O}(\frac{1}{z^{30}})
 \eean
and then using below the trivial Wronskian relation
$-a^2b'=\{a,ab\}$, one computes
 \bean
 \lefteqn{T_2T_5}
  \\&=&
   \frac{1}{z^{31}} T_2^{0} T_5^{0}
    +O\left(\frac{1}{z^{30}}\right)
  \\
  &=&
  \frac{1}{32^3(s_2\!-\!s_1)^4(s_m\!-\!s_{m-1})^3z^{31}}
  \left( -\frac{3}{32}(\XR_{-1}^2{\cal T}_{-1}\BQ_z)^2
  (\XR_{-1}^3\BQ_z)
  \right)
  +
{O}(\frac{1}{z^{30}})
 \\
  &=&
  \frac{1}{32^3(s_2\!-\!s_1)^4(s_m\!-\!s_{m-1})^3z^{31}}
  \left\{ \XR_{-1}^2{\cal T}_{-1}\BQ_z ,
    \frac{3}{32}(\XR_{-1}^2 {\cal T}_{-1}\BQ_z)
  (\XR_{-1}^2\BQ_z)
  \right\}_{\XR_{-1}}
 \\&& +
{O}(\frac{1}{z^{30}}). \eean
Adding these two contributions, one finds 
%
 \bea
   0&=&-T_1\cdot T_2+T_3\cdot T_4\no\\
   &=&\{T_2,T_3\}_{{C}_1}+T_2T_5
 \no\\
 &=&\frac{1}{32^3(s_2 - s_1)^4(s_m- s_{m-1})^3z^{31}}
  \no\\
 &&
 \Biggl\{\XR^2_{-1}{\cal T}_{-1}\BQ_z,~
  \frac{1}{16}\{
\XR_{-1}{\cal T}_{-1} \BQ_z,\XR^2_{-1}
\BQ_z\}_{\XR_{-1}}
 \no\\
& & \hspace*{ 0.5cm}
 +\left({\cal T}_0-1+
\frac{\XR_0}{2}\right)\XR_{-1}^2 \BQ_z-
2(\tilde\XR_{-1}\XR_{-1}-{\cal T}^2_{-1}){\cal T}_{-1}
\BQ_z
      \Biggr\}_{ \XR_{-1}}+
{O}(\frac{1}{z^{30}})\no\\
 \label{final}\eea
  In \cite{TW}, Tracy and Widom show that the extended
  kernel for the non-intersecting Brownian motions tend
  to the extended Pearcey kernel uniformly in each
  bounded interval. Since $\BQ_z$ is the log of the
  Fredholm determinant for that kernel, it follows that
  $$\lim_{z\rightarrow 0}\BQ_z=\BQ_0.$$
   Then taking the limit when $z\rightarrow 0$ in (\ref{final}) leads to
 the PDE for $\BQ_0$, ending the proof of Theorem 0.2.
 \qed


\newpage

\newpage

\section{Appendix: evaluation of the integral over the full range}

    \begin{lemma}

        The following integral can be evaluated
        explicitly:
%
%
\bea 
%
\tau_{k_1k_2}(\BR)&= & \frac{1}{k_1!k_2!} \int_{
(\BR^m)^{k_1+k_2}} \Dt_{k_1+k_2}(y^{(1)})
   \prod_{i=1}^{m}dy_i^{(\ell)}
 \no\\
   && \Dt_{k_1}(y^{(m)'})
\prod_{i=1}^{k_1}e^{-\frac{1}{2}\sum_{\ell=1}^{m}
y_i^{(\ell)2}
  + \sum_{\ell=1}^{m-1}c_{\ell}y_i^{(\ell)}y_i^{(\ell+1)}
   +\alpha y_i^{(m)}  }\no\\
 &&   \Dt_{k_2}(y^{(m)''})
  \!\!\prod_{i=k_1+1}^{k_1+k_2}\!\!
   e^{-\frac{1}{2}\sum_{\ell=1}^{m}
y_i^{(\ell)2}
  + \sum_{\ell=1}^{m-1}c_{\ell}y_i^{(\ell)}y_i^{(\ell+1)}
   -\alpha y_i^{(m)}}
%
  \no\\
 &=&
 \det\left(\begin{array}{c}
\left(\displaystyle{\mu_{i,j}(\al)}\right)_{\tiny{\begin{array}{l}
      0\leq i\leq k_1-1\\
      0\leq j\leq n-1
      \end{array}}}\\
\\
\left(\displaystyle{\mu_{i,j}(-\al)}\right)_{\tiny{\begin{array}{l}
      0\leq i\leq k_2-1\\
      0\leq j\leq n-1
      \end{array}}}
\end{array}\right)
  \no\\ & =&
   c_{k_1k_2} \alpha ^{k_1k_2} 
    \left({\prod_1^{m-1}c_i}\right) ^{-\frac{k_1+k_2}{2}}
  e^{ -\frac{ \alpha^2 }{2}(k_1+k_2)J_{mm}}
 ( J_{1m})^{\frac{1}{2} (k_1+k_2)^2}
  \no\\ \eea
where \bean
 \mu_{ij}(\al)&=&\int_{\BR_m}x_1^jx_m^i
  e^{-\frac{1}{2}
  \left(\sum_1^mx_j^2-2\sum_1^{m-1}c_jx_jx_{j+1}\right)+\alpha
  x_m} dx_1\ldots dx_m
  \\
 c_{k_1k_2}&=&(-2 )^{k_1k_2}
(2\pi)^{\frac{m}{2}(k_1+k_2)} \prod_0^{k_1-1}j!
  \prod_0^{k_2-1}j! ~.
  \eean
Therefore
\bean \left( \bar{\cal A}_2\bar{\cal C}_1+
J_{1m}\frac{\pl}{\pl \alpha}\right)
  \log \tau_{k_1k_2}(\BR)
  &=&J_{1m}\frac{\pl}{\pl \alpha}\left( k_1k_2\log \alpha
-\frac{\alpha^2}{2}(k_1+k_2)J_{mm}\right)
 \\
  &=&J_{1m}\left(\frac{k_1k_2}{\al}-\al (k_1+k_2)
  J_{mm}\right),
\eean
 \bean - \left( \bar{\cal C}_2 + 2\alpha
J_{1m}\bar{\cal C}_1\right){\cal A}_1 \log
\tau_{k_1k_2}(\BR)
 &=&  \frac{1}{2}c_1 \frac{\pl}{\pl c_1}\frac{\pl}{\pl \alpha}
 \left( k_1k_2\log \alpha
-\frac{\alpha^2}{2}(k_1+k_2)J_{mm}\right)
 \\
 &=&\alpha(k_1+k_2)J_{1m}^2.
 \eean
        \end{lemma}

\proof From the explicit evaluation of the zero
moment\footnote{Using
$$
 \int_{\BR^m}e^{-\frac{1}{2}\la Qx,x\ra+\la \ell,x\ra}
 dx_1\ldots dx_m
 =
  \frac{(2\pi)^{m/2}}{\sqrt{\det Q}}
e^{\frac{1}{2}\la Q^{-1}\ell,\ell\ra}
 ,$$
 for $Q:=-J^{-1}$ and $\ell:=(\gamma,0,\ldots,0,\al)$.}
\bean \mu_{00}(  \al,\gamma)&:=& \int_{\BR^m}
 e^{-\frac{1}{2} \left(\sum_1^m
 x_j^2-2\sum_1^{m-1}c_jx_jx_{j+1}\right)+\gamma x_1+\al
 x_m} dx_1\ldots dx_m
 \\&=&
  \frac{(2\pi)^{m/2}}{\sqrt{\det J^{-1}(c_1,\ldots,c_{m-1})}}
e^{-\frac{1}{2}(J_{mm}\al^2+2\al \gamma
J_{1m}+J_{11}\gamma^2)}
 ,\eean
one deduces the other moments by derivation,
 \bean
 \mu_{ij}(\pm \al)&=&
  \left(\pm \frac{\pl}{\pl \al}\right)^i
   \left( \frac{\pl}{\pl \gamma}\right)^j
   \mu_{00}(\pm \al,\gamma)\Bigr|_{\gamma=0}
   \\
   &=&\mu_{00}(\pm \al,\gamma)
   (\pm 1)^{i} A^i B_{\pm}^j(1)\Bigr|_{\gamma=0}
   \\
   &=&\mu_{00}(\pm \al,\gamma)
   (\pm 1)^{i} A^i p_j(\pm \al)
   ,\eean
where
 \bean A&:=&\mu_{00}(\pm \al,\gamma)^{-1}\frac{\pl}{\pl \al}
  \mu_{00}(\pm \al,\gamma)\Bigr|_{\gamma=0}
 \\
 &=&\frac{\pl}{\pl \al}-J_{mm}\al\\ \\
 B_{\pm}&:=&
  \mu_{00}(\pm \al,\gamma)^{-1}\frac{\pl}{\pl \gamma}
  \mu_{00}(\pm \al,\gamma)
  \\& =&
  \frac{\pl}{\pl \gamma}-J_{11}\gamma\mp\al J_{1m}
 \\
 p_j(\al)&:=& B_+^j(1)\Bigr|_{\gamma=0}
 .\eean
 The following holds:
$$p_{2i}(\al)=\mbox{even polynomial,}\qquad p_{2i+1}(\al)=\mbox{odd
polynomial of $\al$},
$$
which is used in equality $\stackrel{\ast\ast}{=}$
below, and
$$A^np_{k}(\al)=p_k^{(n)}+\beta_1(\al)p_k^{(n-1)}+\beta_2(\al)p_k^{(n-2)}+\ldots
+\beta_np_k,
$$
where $p_k^{(n)}:=\left(\frac{d}{d\alpha}\right)^np_k$
and where $\beta_i(\al)$ are polynomials in $\al$,
independent of $k$; this feature is used in equality
$\stackrel{\ast}{=}$ below. The equality
$\stackrel{\ast\ast\ast}{=}$ hinges on the identity
$$
\det\left(\begin{array}{c} \left(
-\al^{j-i}\right)_{\tiny{\begin{array}{l}
      1\leq i\leq k_1\\
      1\leq j\leq n
      \end{array}}}\\
\\
\left( \al^{j-i}\right)_{\tiny{\begin{array}{l}
       1\leq i\leq k_2\\
      1\leq j\leq n
      \end{array}}}
\end{array}\right)=c''_{k_1k_2}\al^{k_1k_2}
$$
for a constant $c''_{k_1k_2}$ depending on $k_1$ and
$k_2$ only.
 Then, setting
$$ \mu:=\mu_{00}(\pm \al,0)=
  \frac{(2\pi)^{m/2}}{\sqrt{\det( -J^{-1}(c_1,\ldots,c_{m-1}))}}
e^{-\frac{1}{2} J_{mm}\al^2 },$$
we compute:
\bean
\lefteqn{\tau_{k_1k_2}(\BR)\big\vert_{t=s=u=\beta =0}}\\
&=&\mu^n~ \det\left(\begin{array}{c}
\left(A^ip_{j}(\al)\right)_{\tiny{\begin{array}{l}
      0\leq i\leq k_1-1\\
      0\leq j\leq n-1
      \end{array}}}\\
\\
\left((-A)^{i }p_{
j}(-\al)\right)_{\tiny{\begin{array}{l}
      0\leq i\leq k_2-1\\
      0\leq j\leq n-1
      \end{array}}}
\end{array}\right)\\
\\
&=&\mu^n~ (-1)^{\frac{k_2(k_2-1)}{2}}
 \det\left(\begin{array}{c}
\left(A^ip_{j}(\al)\right)_{\tiny{\begin{array}{l}
      0\leq i\leq k_1-1\\
      0\leq j\leq n-1
      \end{array}}}\\
\\
\left(A^ip_{j}(-\al)\right)_{\tiny{\begin{array}{l}
      0\leq i\leq k_2-1\\
      0\leq j\leq n-1
      \end{array}}}
\end{array}\right)\\
\\
&\stackrel{\ast}{=}& \mu^n~
(-1)^{\frac{k_2(k_2-1)}{2}}\det
 \left(\begin{array}{c}
 \left(p_{j}^{(i)}(\al)\right)_{\tiny{\begin{array}{l}
      0\leq i\leq k_1-1\\
      0\leq j\leq n-1
      \end{array}}}\\
\\
\left( p_{j}^{(i)}(-\al)\right)_{\tiny{\begin{array}{l}
      0\leq i\leq k_2-1\\
      0\leq j\leq n-1
      \end{array}}}
\end{array}\right)\\
\\
&\stackrel{\ast\ast}{=}& c_{k_1k_2}\mu^n \det\left(
 \begin{array}{c}
 \left(((-J_{1m}\al)^{j-1})^{(i)}\right)_{\tiny{\begin{array}{l}
      0\leq i\leq k_1-1\\
      1\leq j\leq n
      \end{array}}}\\
\\
\left(((J_{1m}\al)^{j-1})^{(i)}\right)_{\tiny{\begin{array}{l}
      0\leq i\leq k_2-1\\
      1\leq j\leq n
      \end{array}}}
\end{array}\right)\\
\\
&=&c_{k_1k_2}\mu^n (J_{1m})^{\frac{n(n-1)}{2}}
\det\left(\begin{array}{c} \left(
-\al^{j-i}\right)_{\tiny{\begin{array}{l}
      1\leq i\leq k_1\\
      1\leq j\leq n
      \end{array}}}\\
\\
\left( \al^{j-i}\right)_{\tiny{\begin{array}{l}
       1\leq i\leq k_2\\
      1\leq j\leq n
      \end{array}}}
\end{array}\right)\\
\\
\eean\bean
&\stackrel{\ast\ast\ast}{=}& c'_{k_1k_2}
   \mu^n (J_{1m})^{\frac{n(n-1)}{2}}\al^{k_1k_2}\\
   &=&c'_{k_1k_2}
   \frac{(2\pi)^{nm/2}}{ {\Dt(c_1,\ldots,c_{m-1})}^{n/2}}
e^{-\frac{n}{2} J_{mm}\al^2 }
(J_{1m})^{\frac{n(n-1)}{2}}\al^{k_1k_2}
\\
 &=& c'_{k_1k_2}(2\pi)^{nm/2}
 \alpha ^{k_1k_2} 
    \left({\prod_1^{m-1}c_i}\right) ^{-\frac{n}{2}}
  e^{-\frac{\alpha^2 }{2}nJ_{mm}}
 ( J_{1m})^{\frac{1}{2} n^2}  .
\eean In order to evaluate the integer $c'_{k_1k_2}$, it
suffices to notice that this constant is independent of
$m$, so that we may choose $m=1$, which was done in
(\cite{AvM4}). This ends the proof of Lemma 8.1.\qed



\begin{thebibliography}{99}


\bibitem{AvM0}
M. Adler, T. Shiota and P. van Moerbeke: {\em Random
matrices, vertex operators and the Virasoro algebra},
Phys. Lett. {\bf A 208}, 67-78, (1995).

\bibitem{AvM1}
M. Adler and P. van Moerbeke: {\sl  The spectrum of
coupled random matrices}, Annals of Math., {\bf 149},
921--976  (1999).

\bibitem{AvM2}
M. Adler and P. van Moerbeke:
 {\sl  Hermitian, symmetric and symplectic random ensembles: PDE's
 for the distribution of the spectrum}, Annals of Math., {\bf
 153}, 149--189 (2001).



\bibitem{AvM3}
M. Adler and P. van Moerbeke: {\sl  PDE's for the joint
distributions of the Dyson, Airy and Sine processes},
The Annals of Probability, {\bf 33}, 1326-1361 (2005).
(arXiv:math.PR/0302329 and math.PR/0403504)


\bibitem{AvM4}
M. Adler and P. van Moerbeke: {\sl PDE's for the
Gaussian ensemble with external source and
  the Pearcey distribution}, Comm. Pure and Appl. Math,
  2006 (arXiv:math.PR/0509047)



\bibitem{AvM5} M. Adler, P. van Moerbeke and P. Vanhaecke:
{\sl  Moment matrices and multicomponent KP, with
applications to random matrix theory}, arXiv:math.PR/
(2006)




\bibitem{AptBleKui}
A. Aptekarev, P. Bleher and A. Kuijlaars: {\em Large $n$
limit of Gaussian random matrices with external source}.
II. {Comm. Math. Phys.}  {\bf 259} (2005) 367--389.,
arXiv: math-ph/0408041.




\bibitem{BleKui1}
P. Bleher and A. Kuijlaars: {\sl  Random matrices with
external source and multiple orthogonal polynomials},
Internat. Math. Research Notices {\bf 3}, 109--129
(2004) (arXiv:math-ph/0307055).


\bibitem{BleKui2}
P. Bleher and A. Kuijlaars: {\sl  Large $n$ limit of
Gaussian random matrices with external source, Part I},
Comm. Math. Phys., {\bf 252}, 43--76  (2004).








\bibitem{Brezin2}
E. Br\'ezin and S. Hikami: {\sl  Correlations of nearby
levels induced by a random potential}, Nuclear Physics
{\bf B 479}, 697--706 (1996).



\bibitem{Brezin3}
E. Br\'ezin and S. Hikami: {\sl  Extension of level
spacing universality}, Phys. Rev., {\bf E 56}, 264--269
(1997).



\bibitem{Brezin4}
E. Br\'ezin and S. Hikami: {\sl  Universal singularity
at the closure of a gap in a random matrix theory},
Phys. Rev., {\bf E 57}, 4140--4149 (1998).


\bibitem{Brezin5}
E. Br\'ezin and S. Hikami: {\sl  Level spacing of random
matrices in an external source}, Phys. Rev., {\bf E 58},
7176--7185 (1998).



 \bibitem{Dyson}
 F.J. Dyson: {\sl  A Brownian-Motion Model for the
 Eigenvalues of a Random Matrix}, Journal of Math. Phys.
 {\bf 3}, 1191--1198 (1962)






\bibitem{Grabiner}
D.J. Grabiner: {\sl  Brownian-Motion in a Weyl chamber,
non-colliding particles and random matrices}, Ann.
Institut H. Poincar\'e {\bf 35}, 177--204 (1999)







\bibitem{Johansson}
K. Johansson: {\sl  Universality of the Local Spacing
distribution in certain ensembles of Hermitian Wigner
Matrices}, Comm. Math. Phys. {\bf 215}, 683--705 (2001)










\bibitem{Karlin}
S. Karlin and J. McGregor: {\sl  Coincidence
probabilities}, Pacific J. Math. {\bf 9}, 1141--1164
(1959).

\bibitem{Okounkov}
A. Okounkov and N. Reshetikhin: {\sl Random skew plane
partitions and the Pearcey process}, math.CO/0503508
(2005).



\bibitem{Pastur}
L.A. Pastur: {\sl  The spectrum of random matrices
(Russian)}, Teoret. Mat. Fiz. {\bf 10}, 102--112 (1972).



\bibitem{Pearcey}
T. Pearcey: {\sl  The structure of an electromagnetic
field in the neighborhood of a cusp of a caustic},
Philos. Mag. {\bf 37}, 311--317 (1946).








\bibitem{TW}
C. Tracy and H. Widom: {\em The Pearcey process}, Comm.
Math. Phys.  {\bf 263}  (2006),  no. 2, 381--400.
arXiv:math. PR /0412005.





































\bibitem{UT} K. Ueno and K. Takasaki: {\em Toda Lattice
Hierarchy}, Adv. Studies in Pure Math. {\bf 4}, 1--95
(1984).


  \bibitem{Zinn}
  P. Zinn-Justin:
 {\sl   Random Hermitian matrices in an external field},  Nuclear Physics {\bf B 497}, 725--732
 (1997).


\bibitem{Zinn1}
P. Zinn-Justin: {\sl   Universality of correlation
functions in Hermitian random matrices in an external
field}, Comm. Math. Phys. {\bf 194}, 631--650 (1998).




\end{thebibliography}
\end{document}